\documentclass[reqno]{amsart}[12pt]
\usepackage{tikz}
\usepackage{pgfplots}
\pgfplotsset{compat=newest}
\usepackage{caption}
\usepackage{float}
\usepackage{mathtools,amsmath,amssymb,amsthm,array,etoolbox}
\usepackage{color}
\usepackage{cases}
\usepackage{ulem}
\usepackage{lineno}
\usepackage{amsmath}
\usepackage{setspace}
 \allowdisplaybreaks
\numberwithin{equation}{section}
\usepackage[left=2.7cm, right=2.7cm, top=2.5cm, bottom=2.5cm]{geometry}
\renewcommand{\baselinestretch}{1.20}

\theoremstyle{definition}

\theoremstyle{remark}

\newcommand{\be}{\begin{equation}}
	\newcommand{\ee}{\end{equation}}
\newcommand{\bse}{\begin{subequations}}
	\newcommand{\ese}{\end{subequations}}
\def\bea{\begin{eqnarray}}
	\def\eea{\end{eqnarray}}
\def\nn{\nonumber}
 \allowdisplaybreaks
\begin{document}
\large
\title [Nonlinear Klein-Gordon System with Quadratic Interaction]{ $L^{2}$-Supercritical Nonlinear Klein-Gordon System with Quadratic Asymmetric Interaction}

\author [X. Dong, H. Duan, Z.H. Gan, J. Wu]{ Xiaojing Dong$^{1}$,~~ Huagui Duan$^{2}$, ~~Zaihui Gan$^{3}$,  Jingyu Wu$^{3}$}
\address{${}^1$ (XD) College of Mathematics and Systems Science, Shandong University of Science and Technology\\ Qingdao, 266590, China.}
\email{dxjmath@sdust.edu.cn(Xiaojing Dong)}
\address{${}^2$ (HD) School of Mathematics, Jilin University \\ Changchun, 130012, Jilin, China.}
 \email{duanhg@jlu.edu.cn(Huagui Duan)}
\address{${}^3$ (ZG,JW) Center for Applied Mathematics, KL-AAGDM, Tianjin University \\ Tianjin, 300072, China.}\thanks{Corresponding author: Zaihui Gan (ganzaihui2008cn@tju.edu.cn)}
\email{ganzaihui2008cn@tju.edu.cn(Zaihui Gan)}
\email{arithwu66@163.com(Jingyu Wu)}

\date{}

\maketitle
\allowdisplaybreaks

\begin{abstract}

In this paper we investigate the global existence, blow-up and standing waves for the $L^{2}$-supercritical nonlinear Klein-Gordon equations with quadratic asymmetric interaction (LSNKG).\\
 \indent First, by introducing a suitable auxiliary functional, using concavity analysis and virial estimates, we obtain a finite time blow-up result for solutions to the Cauchy problem of (LSNKG) when the initial energy is negative. Next, by defining appropriate functionals, manifolds and a constrained variational problem, we employ variational method and the Lagrange multiplier method to derive the existence of ground state solutions for the corresponding  nonlinear elliptic 
 system, thereby to establish the existence of standing wave with the ground state for (LSNKG). Then, using the variational characterization of the ground state solutions and constructing invariant sets under the flow generated by the Cauchy problem for (LSNKG), we combine the potential well argument with concavity analysis to establish a sharp threshold between blow-up in finite time and global existence. Finally, by exploiting the variational characterization  of the ground state, introducing appropriate scalings, and choosing suitable initial data based on the ground state, we justify the instability of standing wave with the ground state for (LSNKG).
 \end{abstract}
{\bf Statements and Declarations:} No conflict of interest exists in the submission of this manuscript. No data was used for the research described in this manuscript.

%
%
%
%
%

\section{Introduction}
 \allowdisplaybreaks
\indent  In this paper, we consider the $L^{2}$-supercritical nonlinear Klein-Gordon equations with quadratic asymmetric interaction:\\
$$
\left\{
\begin{array}{lll}
		\phi_{tt} - \Delta \phi + m_1^2 \phi = 2 a_1 m_1 \phi \psi,&t > 0, \; x \in \mathbb{R}^{N}, \quad(1.1a)  \\[0.3cm]
		\psi_{tt} - \Delta \psi + m_2^2 \psi = 2 a_2 m_2 \phi^2, &t > 0, \; x \in \mathbb{R}^{N} , \quad(1.1b)
\end{array}
\right.\eqno(1.1)
   $$
\\
with initial data\\
$$
\left\{
\begin{array}{lll}
		\phi(0,x)=\phi_0(x), \; \phi_t(0,x)=\phi_1(x),\; x\in\mathbb{R^N},\\\\
		\psi(0,x)=\psi_0(x), \; \psi_t(0,x)=\psi_1(x),\; x\in\mathbb{R^N}.
	\end{array}
\right.\eqno(1.2)
   $$
   \\
Here, $2\leq N\leq 5$, $(\phi,\psi)=(\phi(t,x),\psi(t,x))$  is a pair of real-valued functions with respect to $(t,x)\in \mathbb{R}^{+}\times\mathbb{R}^{N}$, $m_1, ~m_2,~ a_1,~ a_2$ are positive real constants. System (1.1) describes the quadratic asymmetric interaction between two scalar fields, and is viewed as a reduced model of the Dirac-Klein-Gordon system concerning proton-proton interactions and the Maxwell-Higgs system which appears in the abelian Higgs model \cite{HIN2011,Miy2021,Ts2003}. In particular, for system (1.1), quadratic asymmetric nonlinear structure poses new challenges to theoretical analysis for the qualitative properties of its solutions.\\ 
\indent In recent years, scholars have gradually turned to the study of the coupled nonlinear wave equations \cite{Ba1998,GGZ2009,GZ2005,Ge1990,HIN2011,
KOS2012,MN2010,ME1988,MR1985,MM1987,NO2001,OTT1995,RR1982,ST1985,
Su2003,TZ1999,Ts2003,Zh2003,ZGG2010}. These researches have formed a complete methodological system for nonlinear wave equations and large amount of them are concerned with symmetric non-quadratic nonlinear interaction. However, in realistic physical scenarios, nonlinear interactions are often non-symmetric and asymmetric potential functions are also widely adopted in models for cosmic phase transitions and domain wall formation. In particular, as is well known, most physical models such as the relativistic superconductor model (\cite{De1964,Sg1985 1}) involve quadratic terms. On the other hand,
compared with symmetric quadratic nonlinear interaction, quadratic asymmetric  terms bring essential mathematical difficulties.
 \\[0.3cm]
 \indent From a mathematical point of view, to the best of our knowledge,
 studies on coupled wave equations involving asymmetric nonlinear interactions, particularly asymmetric quadratic interactions, remain limited so far.  This paper is devoted to the studies involving finite time blow-up, existence and instability of standing waves with ground state as well as the sharp threshold of global existence for system (1.1) in the $L^{2}-$supercritical case.
 \\[0.3cm]
\indent Let us now recall some known results for system (1.1). When $(\phi,\psi)$ is a $\mathbb{C}^{2}$-valued unknown function with respect to $(t,x)$, let $ a_1 m_1=2 a_2 m_2=1 $, system (1.1) takes the following form\\
$$
\left\{
\begin{array}{lll}
		\phi_{tt} - \Delta \phi + m_1^2  \phi = 2 \bar{\phi} \psi,&t > 0, \; x \in \mathbb{R}^{N},  \\[0.3cm]
		\psi_{tt} - \Delta \psi + m_2^2 \psi =   \phi^2, &t > 0, \; x \in \mathbb{R}^{N}.
\end{array}
\right.\eqno(1.1^*)
   $$
\\
$(1.1^*)$ has the standing wave solutions of the form:
$$(\phi(t,x),\psi(t,x))=(e^{i\omega t}u(x),e^{2i\omega t}v(x)),\eqno(1.2^*)$$
where $(u,v)\in H^1(\mathbb{R^N})\times H^1(\mathbb{R^N})$ is a $\mathbb{C}^{2}$-valued function, which satisfies the elliptic equations below:\\
$$
\left\{
\begin{array}{lll}
		-\Delta u+(m_1^2-\omega^2)u=2 \bar{u}v,&\;
\\[0.3cm]
		-\Delta v+(m_2^2-4\omega^2)v= u^2.&\;
	\end{array}
\right.\eqno(1.3^*)
   $$
   \\
   As $ (u,v)\in H^1(\mathbb{R^N})\times H^1(\mathbb{R^N})$ is a $\mathbb{R}^{2}$-valued function, Miyazaki \cite{Miy2021} adopted the argument in \cite{OH2007} for the single nonlinear Klein-Gordon equation to establish the strong instability of the standing wave solution for $(1.1^*)$. On one hand, for $N=4,5$ ($L^2$-critical and $L^2$-supercritical cases), Miyazaki \cite{Miy2021} proved that the standing wave solution $(e^{i\omega t}u(x),e^{2i\omega t}v(x))$ of $(1.1^*)$ is strongly unstable for all $\omega\in \mathbb{R}$ satisfying $\omega^2<\min\left\{m_1^2,\frac{m_2^2}{4}\right\}$. On the other hand, for $N=2,3$ ($L^2$-subcritical case) and $m_2=2m_1$ (mass resonance), he proved in \cite{Miy2021} that the standing wave solutions $(e^{i\omega t}u(x),e^{2i\omega t}v(x))$ of $(1.1^*)$ possess strong instability under the condition $$(5-N)\omega^2<m_1^2 .\eqno(1.3a)$$
   However, for $N=2,3$ ($L^2$-subcritical case), without requiring the mass resonance condition $m_2=2m_1$, drawing on the method proposed in the study of the nonlinear Schr\"{o}dinger equations with harmonic potential \cite{OH2018}, Ohta in \cite{OH2025} established the strong instability of the standing wave solution $(e^{i\omega t}u(x),e^{2i\omega t}v(x))$ of $(1.1^*)$ with $\omega\in \mathbb{R}$ and
   $$(5-N)\omega^2<\min\left\{m_1^2,\frac{m_2^2}{4}\right\}.\eqno(1.3b)$$
As Ohta stated in \cite{OH2025}, drawing on the approach adopted by Shatah in \cite{Sg1983} for investigating the single Klein-Gordon equation, condition (1.3a) is optimal under the mass resonance condition $m_2=2m_1$.
\\[0.3cm]
\indent It should be mentioned that some works for the single Klein-Gordon equation and for the coupled nonlinear wave equations to elaborate on the research background and significance of the problem investigated in this paper. As fundamental works of the single quadratic nonlinear Klein-Gordon equation,
 Shatah \cite{Sg1985 1} studied the asymptotic behavior by constructing
a transformation and established global existence. On the other hand,
 shatah \cite{Sg1983,Sg1985 2,Sg1985 3} proved the stability and instability of the standing wave solutions. In the research of coupled wave systems, there have been many researches in recent years \cite{Ba1998,GGZ2009,GZ2005,Ge1990,HIN2011,
 KOS2012,MN2010,ME1988,MR1985,MM1987,NO2001,OTT1995,RR1982,ST1985,
Su2003,TZ1999,Ts2003,Zh2003,ZGG2010}. Medeiros-Menzala \cite{ME1988} and Miranda-Medeiros \cite{MM1987} adopted the Galerkin method and energy estimates to prove the existence and uniqueness of weak solutions for coupled nonlinear Klein-Gordon equations describing interactions between charged mesons and electromagnetic fields. Nakamura-Ozawa \cite{NO2001} applied Strichartz estimates on Besov spaces to study nonlinear terms with exponential growth and established scattering theory for certain exponents. Sunagawa \cite{Su2003} studied global small-amplitude solutions for cubic coupled systems under non-resonance conditions. Katayama-Ozawa-Sunagawa \cite{KOS2012} extended relevant results to two-dimensional quadratic systems by normal form transformations and null conditions.
 Bachelot \cite{Ba1998} analyzed the connection between Lorentz invariance and null conditions for the Dirac-Klein-Gordon system. Hayashi-Ikeda-Naumkin \cite{HIN2011} proved the existence of wave operators with low regularity. Tsutsumi \cite{Ts2003} introduced normal forms for the Maxwell-Higgs system. Zhang \cite{Zh2003} defined action functionals and Nehari manifolds to study standing waves for coupled Klein-Gordon equations. Gan-Guo-Zhang \cite{GGZ2009,GZ2005,ZGG2010} extended relevant methods to the Klein-Gordon-Zakharov system and other coupled wave models. Ozawa-Tsutaya-Tsutsumi \cite{OTT1995} transformed quadratic nonlinearities into cubic terms via normal forms for the Klein-Gordon-Zakharov system. Masmoudi-Nakanishi \cite{MN2010} investigated the non-relativistic limit of the system and discovered radiation damping caused by nonlinear resonance.\\

    \indent In the present paper, without requiring mass resonance restriction $m_2=2m_1$, we systematically investigate the qualitative properties of solutions to equations (1.1) for the $L^2$-supcritical case $N=5$, focusing not only on the instability of standing wave solution but also on other qualitative properties of solutions including finite time blow-up, the existence of standing wave with ground state, the sharp conditions for blow-up and global existence for the Cauchy problem of equations (1.1). Although these topics have been widely investigated for the nonlinear Klein-Gordon equations with symmetric coupling and for the single case, they  are novel to some extent for (1.1) as far as we know.\\
\indent Based on these works \cite{GV1985,GV1989,Miy2021,OH2025}, by virtue of the Strichartz estimate \cite{PH1984,WB1998} and the classical semigroup theory \cite{Se1963Nsg}, the local well-posedness for the Cauchy problem (1.1)-(1.2) in the energy space $H^1(\mathbb{R^N})\times H^1(\mathbb{R^N})\times L^2(\mathbb{R^N})\times L^2(\mathbb{R^N})$ with $2\leq N\leq 5$ is given as follows. \\
\\
\textbf{Proposition 1.1}. {\sl \quad Suppose that $m_1>0, m_2>0, a_1>0, a_2>0$ and $2\leq N\leq 5$. For any $
	(\phi_0,\psi_0)\in H^1(\mathbb{R^N})\times H^1(\mathbb{R^N}), \;(\phi_1,\psi_1)\in L^{2}(\mathbb{R^N})\times L^{2}(\mathbb{R^N})
 $, there exists a maximal existence time $T\in (0,+\infty)$ such that the Cauchy problem (1.1)-(1.2) admits a unique solution $(\phi,\psi)$ satisfying\\
$$
	(\phi,\psi)\in C\left([0,T);\;H^1(\mathbb{R^N})\times H^1(\mathbb{R^N})\right),~~(\phi_{t},\psi_{t})\in C\left([0,T);\; L^{2}(\mathbb{R^N})\times  L^{2}(\mathbb{R^N})\right),
$$
\\
and either $T=\infty$ (global existence), or  $T<\infty$ and then\\
$$
	\lim_{t\to T}\left(\|\phi\|_{H^1(\mathbb{R^N})}
+\|\psi\|_{H^1(\mathbb{R^N})}+\|\phi_{t}\|_{L^{2}(\mathbb{R^N})}
+\|\psi_{t}\|_{L^{2}(\mathbb{R^N})}\right)=+\infty~\mbox{(finite~time~blow-up)}.\eqno \Box
$$
}
 \\
\indent We now define the energy functional of the Klein-Gordon equations (1.1) as\\
\\
 $$
\left.
	\begin{array}{ll}
		E(\phi,\psi,\phi_t,\psi_t)
		&\displaystyle=a_2m_2\int_{\mathbb{R^N}}\phi_t^2dx+\frac{a_1m_1}{2}
\int_{\mathbb{R^N}}\psi_t^2dx+a_2m_2\int_{\mathbb{R^N}}|\nabla\phi|^2dx
\\[0.5cm]	&\displaystyle\quad+\frac{a_1m_1}{2}\int_{\mathbb{R^N}}|\nabla\psi|^2dx+a_2m_1^2m_2
\int_{\mathbb{R^N}}|\phi|^2dx\\[0.5cm] &\displaystyle\quad+\frac{a_1m_2^2m_1}{2}\int_{\mathbb{R^N}}|\psi|^2dx-2a_1a_2m_1m_2\int_{\mathbb{R^N}}\phi^2\psi dx.
\end{array}
\right.\eqno (1.3)
$$
\\[0.3cm]
By a multiplier technique, (1.1) admits the following conservation laws in the energy space $H^1(\mathbb{R^N})\times H^1(\mathbb{R^N})\times L^2(\mathbb{R^N})\times L^2(\mathbb{R^N})$:\\
 \\
\textbf{Lemma 1.1}.{\sl \quad Let $(\phi,\psi,\phi_t,\psi_t)$
 be a smooth solution to the Cauchy problem (1.1)-(1.2). Then the energy $E(\phi,\psi,\phi_t,\psi_t)$ satisfies the following conservation:\\
$$
		E(\phi,\psi,\phi_t,\psi_t)
		 =E(\phi_0,\psi_0,\phi_1,\psi_1)\triangleq E_{0}.\eqno(1.4)$$
}\\
{\bf Remark 1.1.}\quad Let
 $$(u,v)_{L^{2}(\mathbb{R}^{N})}=\mbox{Re} \int_{\mathbb{R}^{N}} u(x)\overline{v(x)}~dx.$$
 \\
 As mentioned in \cite{Miy2021,OH2025}, if $(\phi,\psi)$ is a $\mathbb{C}^{2}$-valued solution of the Cauchy problem (1.1)-(1.2) in the sense of Proposition 1.1, then it also satisfies the conservation law of charge:\\
$$Q^*(\phi,\psi,\phi_t,\psi_t)=Q^*(\phi_{0},\psi_{0},\phi_1,\psi_1),$$
\\
where $Q^*(u_{1},u_{2},v_{1},v_{2})$ is defined by\\
$$
\left.
\begin{array}{lll}
		Q^*(u_{1},u_{2},v_{1},v_{2})&=\displaystyle
(v_{1},iu_{1})_{L^{2}(\mathbb{R}^{N})}+2(v_{2},iu_{2})_{L^{2}(\mathbb{R}^{N})}
\\[0.3cm]
		 &=\displaystyle \mbox{Im}\int_{\mathbb{R^N}}\overline{u_{1}}v_{1}dx+
2\mbox{Im}\int_{\mathbb{R^N}}\overline{u_{2}}v_{2}dx.
\end{array}
\right.
$$
 \\
However, for our case, $(\phi,\psi)$ is a $\mathbb{R}^{2}$-valued solution of the Cauchy problem (1.1)-(1.2), it yields that the conserved charge is zero.\hfill$\Box$
\\[0.3cm]
\indent The outline of the paper is as follows. In Section 2 , we give some preparatory materials and state the main results. In Section 3, we investigate finite time blow-up solutions to the Cauchy problem (1.1)-(1.2) (the proof of Theorem 2.1). In Section 4, we solve a constrained variational problem and apply the variational method to establish the existence of standing wave solutions with ground state for system (1.1) (the proof of Theorem 2.2). In Section 5, combining the potential well argument and the concavity analysis, we derive a sharp condition for the finite time blow-up and global existence of solutions to the Cauchy problem (1.1)-(1.2) (the proof of Theorem 2.3). In the final section, we employ the variational properties of ground state solutions and adopt appropriate scaling transformations to prove the instability of standing wave with ground state for system (1.1) (the proof of Theorem 2.4).\\

\section{Preliminaries and Main results}
\indent~~In this section, we give some preparatory materials and state the main results.
\subsection{Preparatory materials}
\quad\\
\indent~~We now first list a compactness lemma in \cite{St1997}.\\
\\
{\bf Lemma 2.1 (Compactness)}. {\sl \quad Let $N\geq 2$ and $u\in H_r^1(\mathbb{R}^N)$. Then
 the injection $H_r^1(\mathbb{R}^N)\hookrightarrow L^p(\mathbb{R}^N)$  is compact, for $2<p<\frac{2N}{N-2}$.}\hfill$\Box$
 \\[0.3cm]
 \indent~~Next, we recall some basic properties of Schwarz symmetrization. We first mention the definition of the Schwarz spherical rearrangement (or symmetrization) of a function (see Berestycki-Lions \cite{BC1981}, Gan-Wang \cite{GW2023}).\\
 \\
{\bf Definition 2.1} (Schwarz symmetrization  \cite{BL1983}).
\\[0.3cm]
 \indent {\sl Let $f\in L^1(\mathbb{R}^N)$ be a nonnegative function, then its Schwarz symmetrized function $f^*$ is the unique spherically symmetric, non-increasing (in $r=|x|$), measurable function such that for all $\alpha>0$,
$m\left\{x\in\mathbb{R}^N: f^*\geq\alpha\right\}=m\left\{x\in\mathbb{R}^N: |f|\geq\alpha\right\}$, where $m$ is the Lebesgue measure.} \hfill$\Box$\\
\\
 \indent We next refer to some properties of the Schwarz symmetrization.\\
 \\
{\bf Proposition 2.1} (Basic properties of Schwarz symmetrization ~\cite{BL1983,GW2023}). \\
\\
\indent {\sl Let $f^*$ and $g^*$ be the Schwarz symmetrization of functions $ |f| $ and $ |g| $, respectively. Then there hold:
\\[0.3cm]
 (1) For every continuous function $F$ such that $F(f)$ is integrable, then $$\int_{\mathbb{R}^N}F(f)dx=\int_{\mathbb{R}^N}F(f^*)dx.$$
\\
 (2)~ Riesz inequality:   Let $f$ and $g$ be in $L^2(\mathbb{R}^N)$, then $$\int_{\mathbb{R}^N}f(x)g(x)dx\leq\int_{\mathbb{R}^N}f^*(x)g^*(x)dx.$$
\\
 (3)~ Let $f$ , $g$ be in $L^2(\mathbb{R}^N)$, then $||f^*-g^*||_{L^2(\mathbb{R}^N)}\leq||f-g||_{L^2(\mathbb{R}^N)}$.
\\[0.3cm]
 (4)~ $\displaystyle\int_{\mathbb{R}^N}|f^*|^pdx=\int_{\mathbb{R}^N}|f|^pdx$ for all~$1\leq p<\infty$ such that $f\in L^p(\mathbb{R}^N)$.
\\[0.3cm]
 (5)~ Let $f$ be in $\mathcal{D}^{1,2}(\mathbb{R}^N)$ if $N\geq 3$~ \Big(respectively, in $H^1(\mathbb{R}^N)$ for any $N$\Big), then $f^*$  belongs to $\mathcal{D}^{1,2}(\mathbb{R}^N)$~ \Big(respectively, to $H^1(\mathbb{R}^N)$\Big), and there holds
 $$\displaystyle\int_{\mathbb{R}^N}|\nabla f^*|^2dx\leq\int_{\mathbb{R}^N}|\nabla f|^2dx .$$
\\[0.3cm]
 (6)~ Let $f_{\lambda}(x)=\lambda^{\frac{N}{2}}f(\lambda x)$,~ then $(f_{\lambda})^*=(f^*)_{\lambda}$.}\hfill$\Box$
 \subsection{Constrained Variational Problem}
\quad\\[0.3cm]
\indent We consider the corresponding steady-state nonlinear system:\\
$$
\left\{
\begin{array}{lll}
		-\Delta u+m_1^2u=2a_1m_1uv,&\;(2.1a)\\\\
		-\Delta v+m_2^2v=2a_2m_2u^2,&\;(2.1b)
	\end{array}
\right.\eqno(2.1)
   $$
   \\
where $(u,v)\in  H^1(\mathbb{R^N})\times H^1(\mathbb{R^N})\setminus\{(0,0)\}$,
 and
\begin{equation}
	(u,v)=(u(x),v(x)):\mathbb{R^N}\times \mathbb{R^N}\to\mathbb{R}\times\mathbb{R}.
	\notag
\end{equation}
Any non-trivial solution $(u,v)\in  H^1(\mathbb{R^N})\times H^1(\mathbb{R^N})\setminus\{(0,0)\}$ of system (2.1) gives a standing wave solution
$$
	(\phi(t,x),\psi(t,x))=(u(x),v(x)),\;t\geq0,x\in\mathbb{R^N}
	\eqno(2.2)
$$
for the original evolution system (1.1).\\
\indent On the other hand, a ground state solution refers to a non-trivial solution with the minimal action among all solutions of the elliptic system. We define the static energy functional $S(u,v)$ and the constraint functional $Q(u,v)$, respectively:\\
 $$
\left.
\begin{array}{lll}
		\displaystyle S(u,v)
		=&\displaystyle a_2m_2\int_{\mathbb{R^N}}|\nabla u|^2dx+\frac{a_1m_1}{2}\int_{\mathbb{R^N}}|\nabla v|^2dx+a_2m_2m_1^2\int_{\mathbb{R^N}}|u|^2dx
\\[0.4cm]	&\displaystyle+\frac{a_1m_1m_2^2}{2}\int_{\mathbb{R^N}}|v|^2dx
-2a_1a_2m_1m_2\int_{\mathbb{R^N}}vu^2dx,
	\end{array}
\right.\qquad\eqno(2.3)
   $$
   \\
$$
\left.
\begin{array}{lll}
		Q(u,v)=&\displaystyle 2a_2m_2\int_{\mathbb{R^N}}|\nabla u|^2dx+a_1m_1\int_{\mathbb{R^N}}|\nabla v|^2dx+2a_2m_2m_1^2\int_{\mathbb{R^N}} u^2dx
\\[0.4cm]
		&\displaystyle +a_1m_1 m_2^2\int_{\mathbb{R^N}}v^2dx-6a_1a_2m_1m_2\int_{\mathbb{R^N}}vu^2dx.
\end{array}
\right.\eqno(2.4)
   $$
 \\[0.3cm]
We further define a manifold $M$ as\\
$$	M=\left\{(u,v)\in H^1(\mathbb{R^N})\times H^1(\mathbb{R^N})\setminus \{(0,0)\}\; : \ Q(u,v)=0\right\}.\eqno(2.5)
$$
\quad\\
We then claim the following.
\\[0.3cm]
{\bf Proposition 2.2}. \quad {\sl The functionals $S(u,v)$ and $Q(u,v)$ are well-defined and belong to the class $C^1$ on $H^1(\mathbb{R^N})\times H^1(\mathbb{R^N})$.}\\[0.3cm]
 \indent We will verify it in the Appendix A.\hfill$\sharp$
 \\[0.3cm]
\indent Let
$$
	G(u,v)=2a_2m_2m_1^2 u  ^2+a_1m_1m_2^2 v ^2-6a_1a_2m_1m_2v u^2.
	\eqno(S^*)
$$
\\
We claim the following preliminary estimates for $G(u,v)$.
 \\[0.3cm]
  {\bf Proposition 2.3}.\quad {\sl Assume  $m_1>0, m_2>0, a_1>0, a_2>0$ and $N=5$. For any solution $(u,v)\in H^1(\mathbb{R^N})\times H^1(\mathbb{R^N})$ to the elliptic system (2.1), there exist $ u^*>0,v^*>0$  such that\\
$$
	G(u^*,v^*)=2a_2m_2m_1^2 (u^*)^2+a_1m_1m_2^2(v^*)^2-6a_1a_2m_1m_2v^*(u^*)^2<0.
	\eqno(S.1)
$$
Here,
$$
	\frac{\partial G(u^*,v^*)}{\partial u^*}=4a_2m_1^2m_2u^*-12a_1a_2m_1m_2v^*u^*,
		\eqno(S.2)
$$
\\
$$
	\frac{\partial G(u^*,v^*)}{\partial v^*}=2a_1m_2^2m_1v^*-6a_1a_2m_1m_2(u^*)^2.
		\eqno(S.3)
$$
}
\\
{ \bf Proof}. Due to $G(0,0)=0$, integrating with respect to $u^*$ on the both
 sides of (S.2) gives\\
$$
\left.
\begin{array}{lll}
			G(u^*,v^*)
			&\displaystyle=\int_{0}^{u^*}(4a_2m_2m_1^2s-12a_1a_2m_1m_2sv^*)ds
    \\[0.4cm]
			&=2a_2m_2m_1^2(u^*)^2-6a_1a_2m_1m_2v^*(u^*)^2+\Phi(v^*),
		\end{array}
\right.\eqno(S.4)
   $$\\
	and hence $\Phi^{'}(v^*)=2a_1m_2^2m_1v^*$,
	that is,
$$
		\Phi(v^*)=\int_{0}^{v^*}2a_1m_2^2m_1sds=a_1m_2^2m_1(v^*)^2.
		\eqno(S.5)
   $$
	 (S.4) and (S.5) then lead to
	$$
		G(u^*,v^*)=2a_2m_2m_1^2(u^*)^2+a_1m_1m_2^2(v^*)^2-6a_1a_2m_1m_2v^*(u^*)^2.
		\eqno(S.6)
   $$
	Since $\{(u,v):G(u,v)=0\}$ is a one-dimensional set, and
	$$
		G(u,u)=2a_2m_2m_1^2u^2+a_1m_1m_2^2u^2-6a_1a_2m_1m_2u^3,
		$$
	for $u_0$ large enough, as $u\geq u_0$, there holds $G(u,u)<0$.\\
	\indent Hence, there exists $u^*>0, v^*>0$ large sufficiently such that
$$ G(u^*,v^*)=2a_2m_2m_1^2(u^*)^2+a_1m_1m_2^2(v^*)^2-6a_1a_2m_1m_2v^*(u^*)^2<0.\eqno \Box		
	$$
\\
\indent In view of Proposition 2.3, we further claim the following conclusion.\\
\\
{\bf Lemma 2.2}.\quad {\sl If $(u_{0},v_{0})\in H^1(\mathbb{R^N})\times H^1(\mathbb{R^N})$ is a solution of (2.1), then\\
$$Q(u_{0},v_{0})=0.$$}\\
{\bf Proof.}\quad Let $u_{\beta}(x)=\beta u_{0}(x),~~v_{\beta}(x)=\beta v_{0}(x)$. Then \\
$$
\left.
\begin{array}{lll}
		\displaystyle S(u_{\beta},v_{\beta})
		&=\displaystyle a_2m_2\int_{\mathbb{R^N}}|\nabla u_{\beta}|^2dx+\frac{a_1m_1}{2}\int_{\mathbb{R^N}}|\nabla v_{\beta}|^2dx
\\[0.4cm]
&~~~~\quad\displaystyle+a_2m_2m_1^2\int_{\mathbb{R^N}}|u_{\beta}|^2dx
 +\frac{a_1m_1m_2^2}{2}\int_{\mathbb{R^N}}|v_{\beta}|^2dx
\\[0.4cm]
&~~~~\quad\displaystyle-2a_1a_2m_1m_2\int_{\mathbb{R^N}}v_{\beta}u_{\beta}^2dx
\\[0.4cm]
&=\displaystyle \beta^{2}\left(a_2m_2\int_{\mathbb{R^N}}|\nabla u_{0}|^2dx+\frac{a_1m_1}{2}\int_{\mathbb{R^N}}|\nabla v_{0}|^2dx\right.
\\[0.4cm]
&~~~~~~~~~\qquad\displaystyle\left.+a_2m_2m_1^2\int_{\mathbb{R^N}}|u_{0}|^2dx
+ \frac{a_1m_1m_2^2}{2}\int_{\mathbb{R^N}}|v_{0}|^2dx\right)
\\[0.4cm]
&~~~~\quad\displaystyle
-2a_1a_2m_1m_2\beta^{3}\int_{\mathbb{R^N}}v_{0}u_{0}^2dx.
	\end{array}
\right.\eqno(B.1)
   $$
\\
Since $(u_{0},v_{0})$ is a solution of (2.1), then $\delta S(u_{0},v_{0})=0$, that ~is,\\
$$\delta_{u_{0}} S(u_{0},v_{0})=0,\quad \delta_{v_{0}} S(u_{0},v_{0})=0,~~ \mbox{or}~~
 \frac{d S(u_{\beta},v_{\beta})}{d \beta}|_{\beta=1}=0,$$\\
 we then obtain
\\
$$
\left.
\begin{array}{lll}
 &\displaystyle 2a_2m_2\int_{\mathbb{R^N}}|\nabla u_{0}|^2dx+ a_1m_1 \int_{\mathbb{R^N}}|\nabla v_{0}|^2dx
 +2a_2m_2m_1^2\int_{\mathbb{R^N}}|u_{0}|^2dx
\\[0.4cm]
&\qquad\qquad\displaystyle
+  a_1 m_1 m_2^2 \int_{\mathbb{R^N}}|v_{0}|^2dx

-6a_1a_2m_1m_2 \int_{\mathbb{R^N}}v_{0}u_{0}^2dx=0.
	\end{array}
\right.\eqno(B.2)
   $$
\\
Then $Q(u_{0},v_{0})=0$.\hfill$\Box$\\
\\
{\bf Lemma 2.3 }.\quad {\sl Let $N=5$. Then $M$ is a $C^{1}-$ hypersurface in $H^1(\mathbb{R^N})\times H^1(\mathbb{R^N})$ bounded away from $(0,0)$.}
\\[0.3cm]
{\bf Proof.}\quad Note that\\
$$	M=\left\{(u,v)\in H^1(\mathbb{R^N})\times H^1(\mathbb{R^N})\setminus \{(0,0)\}\; , Q(u,v)=0\right\},\eqno(B.3)
$$
\\
from Proposition 2.2 and Proposition 2.3, it follows that $Q(u,v)$ is a $C^{1}$ functional on $H^1(\mathbb{R^N})\times H^1(\mathbb{R^N})$. We then claim:
\\[0.3cm]
\indent {\bf Conclusion I:}\quad {\sl For $(u_{0},v_{0})\in M$,
$$\left(\delta_{u_{0}} Q(u_{0},v_{0}),\delta_{v_{0}} Q(u_{0},v_{0})\right)\neq (0,0).\eqno(B.4)$$}

{In fact, consider $(u_{0},v_{0})\in H^1(\mathbb{R^N})\times H^1(\mathbb{R^N})$
such that \\
$$\delta_{u_{0}} Q(u_{0},v_{0})=0,~~\delta_{v_{0}} Q(u_{0},v_{0})= 0.$$
\\
This implies that $(u_{0},v_{0})$ satisfies the following system
\\
$$
\left\{
\begin{array}{lll}
		-4a_2m_2\Delta u_{0}+4a_2m_2m_1^2u_{0}-12a_1a_2m_1m_2u_{0}v_{0}=0,
\\[0.4cm]
		-2a_1m_1\Delta v_{0}+2a_1m_1m_2^2v_{0}-6a_1a_2m_1m_2u_{0}^2=0.
	\end{array}
\right.\eqno(B.5)
   $$
   \\
   That is,\\
 $$
\left\{
\begin{array}{lll}
		- \Delta u_{0}+ m_1^2u_{0}=3a_1 m_1 u_{0}v_{0},
\\[0.4cm]
		- \Delta v_{0}+ m_2^2v_{0}=3 a_2m_2u_{0}^2.
	\end{array}
\right.\eqno(B.6)
   $$
   \\
By a similar argument as in the definition of $S(u,v)$ (see (2.3)), we formulate functional $S^*(u_{0},v_{0})$ as\\
 $$
\left.
\begin{array}{lll}
		\displaystyle S^*(u_{0},v_{0})
		=&\displaystyle a_2m_2\int_{\mathbb{R^N}}|\nabla u_{0}|^2dx+\frac{a_1m_1}{2}\int_{\mathbb{R^N}}|\nabla v_{0}|^2dx
\\[0.4cm]
&\displaystyle+a_2m_2m_1^2\int_{\mathbb{R^N}}|u_{0}|^2dx
+\frac{a_1m_1m_2^2}{2}\int_{\mathbb{R^N}}|v_{0}|^2dx
\\[0.4cm]
&\displaystyle-3a_1a_2m_1m_2\int_{\mathbb{R^N}}v_{0}u_{0}^2dx.
	\end{array}
\right.\eqno(B.7)
   $$
   \\
Then by Lemma 2.2 applied to (B.6) and (B.7), we attain\\
$$
\left.
\begin{array}{lll}
		 		 &\displaystyle 2a_2m_2\int_{\mathbb{R^N}}|\nabla u_{0}|^2dx+ a_1m_1 \int_{\mathbb{R^N}}|\nabla v_{0}|^2dx
\\[0.4cm]
&\displaystyle\qquad+2a_2m_2m_1^2\int_{\mathbb{R^N}}|u_{0}|^2dx
+ a_1m_1m_2^2 \int_{\mathbb{R^N}}|v_{0}|^2dx
\\[0.4cm]
&\displaystyle\qquad-9a_1a_2m_1m_2\int_{\mathbb{R^N}}v_{0}u_{0}^2dx=0.
	\end{array}
\right.\eqno(B.8)
   $$
   \\
Then if $Q(u_{0},v_{0})=0$, it yields by (B.8) that
$$
\left.
\begin{array}{lll}
		 		 &\displaystyle 2a_2m_2\int_{\mathbb{R^N}}|\nabla u_{0}|^2dx+ a_1m_1 \int_{\mathbb{R^N}}|\nabla v_{0}|^2dx
\\[0.4cm]
&\displaystyle\qquad+2a_2m_2m_1^2\int_{\mathbb{R^N}}|u_{0}|^2dx
+ a_1m_1m_2^2 \int_{\mathbb{R^N}}|v_{0}|^2dx=0,
	\end{array}
\right.\eqno(B.9)
   $$
   \\
and this gives $\|u_{0}\|_{H^{1}(\mathbb{R^N})}+\|v_{0}\|_{H^{1}(\mathbb{R^N})}=0$, which is contradictory to $(u_{0},v_{0})\in M$. Hence {\bf Conclusion I} holds, and $M$ is a $C^{1}$-hypersurface.}

\medskip

On the other hand, from Young's inequality and Gagliardo-Nirenberg inequality, it follows that \\
$$
\left.
\begin{array}{lll}
		 \displaystyle\int_{\mathbb{R^N}}v u ^2 dx
		&\displaystyle\leq C\left(\| u \|_{L^p(\mathbb{R^N})}^p+\|v \|_{L^{2q}(\mathbb{R^N})}^{2q}\right)
\\[0.5cm]
		&\displaystyle \leq C\|u \|_{L^2(\mathbb{R^N})}^{p-\frac{N}{2}(p-2)}\|\nabla u \|_{L^2(\mathbb{R^N})}
^{\frac{N}{2}(p-2)}
\\[0.5cm]
&\displaystyle\quad+C\|v \|_{L^2(\mathbb{R^N})}^{2q-\frac{N}{2}(2q-2)}\|\nabla v \|_{L^2(\mathbb{R^N})}^{\frac{N}{2}(2q-2)},
  \end{array}
\right.\eqno(B.10)
$$
\\
where for $N=5$,\\
$$
\left\{
\begin{array}{lll}
	\displaystyle\frac{1}{p}+\frac{1}{q}=1~(p>2,~q>2),~	\frac{N}{4}(p-2)>1,~~ \frac{4N}{N+2}<2q<\frac{4N+8}{N+4},\\[0.5cm]		\displaystyle\frac{N}{4}(2q-2)>\frac{N}{4}(\frac{4N}{N+2}-2)=\frac{N}{4}(\frac{2N-4}{N+2})=\frac{N}{2}(\frac{N-2}{N+2})>1.\\
\end{array}
\right.\eqno(B.11)
   $$
   \\
Let $\theta=\max\left\{\frac{N}{4}(p-2),\frac{N}{4}(2q-2)\right\}>1$ and let
\\
$$
\left.
\begin{array}{lll}
		&\displaystyle A(u ,v )= 2a_2m_2\int_{\mathbb{R^N}}|\nabla u |^2dx+a_1m_1\int_{\mathbb{R^N}}|\nabla v |^2dx
\\[0.4cm]
&\qquad\qquad \displaystyle+2a_2m_1^2m_2\int_{\mathbb{R^N}} u ^2dx
		 +a_1m_2^2m_1\int_{\mathbb{R^N}}v ^2dx.
\end{array}
\right.\eqno(B.12)
   $$
   \\
Then in view of Sobolev embedding, by (B.10) and (B.12) there holds\\
$$
Q(u ,v )\geq A-c_{0}A^{\alpha}\quad for \quad \alpha>1.$$\\
Now for $0<\|u\|_{H^1(\mathbb{R^N})}+\|v\|_{H^1(\mathbb{R^N})}<\varepsilon$ as $\varepsilon$ small, there holds $Q(u,v)>0$, and this implies that $M$ is bounded away from $(0,0)$.\hfill$\Box$\\
\\
\indent Based on Proposition 2.2 and Proposition 2.3, by Lemma 2.2 and Lemma 2.3, we obtain
\\[0.3cm]
{\bf Proposition 2.4}. {\sl $M$ is non-empty.}\hfill$\Box$\\[0.3cm]
\indent We then
 formulate the constrained minimization problem:\\
$$
	m=\inf_{(u,v)\in M}S(u,v).
	\eqno(2.6)
$$
 \subsection{Statement of Main Results}
  \quad\\
  \\
  \indent Before we state the main results of this paper, we formulate several pivotal definitions for the system (1.1) and system (2.1). Firstly, we give a difinition of the solutions to (2.1).\\
  \\
  {\bf Definition 2.2} (Solutions to (2.1) \cite{HOT2013,Miy2021,OH2025}). {\sl We call that a pair of real-valued functions $(u,v)\in  H^1(\mathbb{R}^N) \times H^1(\mathbb{R}^N)$ is a solution to (2.1) if\\
$$
\left\{
\begin{array}{lll}
&\displaystyle\int_{\mathbb{R}^N} \nabla u\cdot \nabla u_1 \,dx
+ m_{1}^2 \int_{\mathbb{R}^N} u u_1 \,dx
= 2a_{1}m_{1}\int_{\mathbb{R}^N} uv u_1 \,dx,
\\[0.5cm]
&\displaystyle\int_{\mathbb{R}^N} \nabla v\cdot \nabla v_1 \,dx
+ m_{2}^2 \int_{\mathbb{R}^N} v v_1 \,dx
= 2a_{2}m_{2}\int_{\mathbb{R}^N} u^2 v_1 \,dx
\end{array}
\right.\eqno(2.6a)
   $$
   \\
for any $(u_{1},v_{1} )\in  C_0^\infty(\mathbb{R}^N) \times  C_0^\infty(\mathbb{R}^N) $. In particular,   $(u,v)$ is a solution to (2.1) if and only if
$
S'\left[(u,v)\right](u_{1},v_{1} ) = 0
$
for any $(u_{1},v_{1})\in H^1(\mathbb{R}^N)\times H^1(\mathbb{R}^N)$.}\hfill$\Box$\\
\\
 \indent Furthermore, we hope to investigate a specified solution to (2.1) named as ground state solution. From the physical viewpoint, an important role is played by the ground state solution of (2.1). The definition of the ground state for (2.1) in the present paper is given below.\\
 \\
  {\bf Definition 2.3} (Ground state solution to (2.1) \cite{HOT2013,Miy2021,OH2025}).  {\sl We call that a pair of real-valued functions $(u,v)\in H^1(\mathbb{R}^N)\times H^1(\mathbb{R}^N)$ is the ground state solution for (2.1) if $(u,v)$ satisfies  \\
  $$S(u,v)\leq S(u^*,v^*)\eqno(2.6b)$$
  \\
  for any $(u^*,v^*)\in H^1(\mathbb{R}^N)\times H^1(\mathbb{R}^N)$
   is a nontrivial critical point of $S(u,v)$. In particular, $(u,v)$ is called a critical point of $S(u,v)$ if $S'\left[(u,v)\right](u_{1},v_{1} ) = 0$ for any $(u_{1},v_{1})\in H^1(\mathbb{R}^N)\times H^1(\mathbb{R}^N)$.}\\
   \indent \quad \hfill$\Box$\\
\\
  {\bf Remark 2.1.}  Note that in \cite{HOT2013}, it is shown that $(1.3^*)$ has positive radially symmetric solutions in $H^1(\mathbb{R}^N)\times H^1(\mathbb{R}^N)$. However, in this paper, we will show that there exist similar solutions for system (2.1).\hfill$\Box$
  \\[0.3cm]
  \indent We then investigate the properties of ground state solutions for (2.1), and utilize them to derive the instability of standing wave solutions for (1.1).
As is well known, the study of the stability and instability for multi-dimensional standing waves is currently regarded as one of the most fascinating research topics in the theory of modern nonlinear waves.\\
\\
 \indent In particular, we give definitions of strong blow-up instability for the standing wave solutions of system (1.1) considered in this paper.\\
 \\
{\bf Definition 2.4} (Strong blow-up instability \cite{Miy2021,OH2025}).  {\sl The standing wave solution $(w(x),\xi(x))$ of (1.1) is called strongly unstable, if for every $\varepsilon>0$, there exists $(\phi_{0},\psi_{0})\in H^1(\mathbb{R^N})\times H^1(\mathbb{R^N})$, and  $(\phi_{1},\psi_{1})=(0,0)$ such that
 \\
$$\|\phi_{0}-w\|_{H^1(\mathbb{R^N})}< \varepsilon,~~\|\psi_{0}-\xi\|_{H^1(\mathbb{R^N})}< \varepsilon,$$
\\
and $T_{max}(\phi_{0},\psi_{0},\phi_{1},\psi_{1})<\infty$, i.e., the solution $(\phi(t),\psi(t))$ of (1.1) with initial data\\
$$(\phi(0),\psi(0),\phi_{t}(0),\psi_{t}(0))
=(\phi_{0},\psi_{0},0,0)$$ blows up in finite time.}\hfill$\Box$\\
\\
{\bf Remark 2.2}.\quad  Generally, from \cite{HOT2013,Miy2021,OH2025} the standing wave solutions for (1.1) is of the form :
\\
$$(\phi(t,x),\psi(t,x))=(e^{i\omega t}u(x),e^{2i\omega t}v(x)),\eqno(2.6b)$$
\\
where $\omega\in \mathbb{R}$ is a constant, and $(u,v)\in H^1(\mathbb{R^N})\times H^1(\mathbb{R^N})$ is a ground state of the stationary system (2.1).  However, in the present paper we only consider the case where $(\phi,\psi)$
 is a pair of real-valued functions in system (1.1) . Consequently, the standing wave solutions in Definition 2.4 are actually a special case of Standing Wave solutions (2.6b), corresponding to $\omega=0$.
 \hfill$\Box$
 \\[0.3cm]
  \indent We now state our main results.
  \\[0.3cm]
   \indent Generally, for a class of sufficiently large data, the solutions of the Cauchy problem for the single nonlinear Klein-Gordon equation blow up in a finite time. Thus in the course of nature, the following problem occurs for the Klein-Gordon system (1.1).
    \\[0.3cm]
{\bf Problem 1:}\quad For the nonlinear Klein-Gordon system  with quadratic asymmetric interaction (1.1), do solutions to its Cauchy problem (1.1)-(1.2) blow up in finite time when the initial data are sufficiently large?\\[0.3cm]
\indent To address this problem, we build upon the local well-posedness theory for the Cauchy problem (1.1)-(1.2) (see Proposition 1.1). Inspired by existing results on finite time blow-up induced by large initial data for the single nonlinear Klein-Gordon equation, we first prove that the solutions to (1.1)-(1.2) blow up in finite time provided that the initial energy is negative. This result verifies  finite time blow-up phenomenon for large initial data and gives an affirmative answer to Problem 1. The detailed content is presented as follows.\\
\\
\textbf{Theorem 2.1} (Finite Time Blow-up).\quad {\sl Assume that $m_1>0,m_2>0,a_1>0,a_2>0$ and $2\leq N\leq 5$. Let $(\phi,\psi,\phi_{t},\psi_{t})\in C([0,T);\;H^1(\mathbb{R^N})\times H^1(\mathbb{R^N})\times L^2(\mathbb{R^N})\times L^2(\mathbb{R^N}))$ be a solution to the Cauchy problem (1.1)-(1.2). If the initial data $(\phi_0,\psi_0,\phi_1,\psi_1)\in H^1(\mathbb{R^N})\times H^1(\mathbb{R^N}))\times L^2(\mathbb{R^N})\times L^2(\mathbb{R^N}))$, and satisfy the following two conditions: \\
\\
(1)\quad $E(\phi_0,\psi_0,\phi_1,\psi_1)=E_{0}<0$;\hfill (2.7)\\
\quad\\
(2)\quad $\displaystyle a_2 m_2\int_{\mathbb{R^N}} \phi_0\phi_1dx
+\frac{a_1m_1}{2}\int_{\mathbb{R^N}}\psi_0\psi_1 dx>0$;\hfill (2.8)\\
\quad\\
then the solution $(\phi,\psi)$ to the Cauchy problem (1.1)-(1.2) blows up in finite time. That is, there exists a time $T\in (0,+\infty)$ such that
$$\lim_{t\to T}\left(\parallel\phi\parallel_{H^1(\mathbb{R^N})}
+\parallel\psi\parallel_{H^1(\mathbb{R^N})}\right)=+\infty.\eqno(2.9)$$}
\\[0.3cm]
\indent In the study of nonlinear wave equations, solutions with special properties can better reflect the inherent characteristics of the equations. For the nonlinear Klein-Gordon system (1.1), standing wave solutions are a special class of particular solutions that are of great physical significance. In particular, standing waves associated with ground states play an essential role in investigating the qualitative properties of solutions to system (1.1).
It is known that a single nonlinear Klein-Gordon equation always admits standing waves with ground state when its nonlinear term satisfies certain conditions. Nevertheless, for the quadratic asymmetric nonlinear Klein-Gordon system (1.1), a natural question arises.
\\[0.3cm]
{\bf Problem 2:}~~ Does system (1.1) possess standing waves with ground state?
\\[0.3cm]
\indent As a matter of fact, the existence of standing waves with ground state for system (1.1) is equivalent to the existence of ground state solutions for its corresponding elliptic system (2.1). Based on formulas (2.3), (2.4), (2.5) and (2.6) in Section 2.2, we can prove the existence of ground state solutions to elliptic system (2.1) by solving the minimization problem (2.6). The following theorem gives an affirmative answer to Problem 2.
\\[0.3cm]
{\bf Theorem 2.2}.\quad  {\sl Suppose $m_1>0,m_2>0,a_1>0,a_2>0$ and $N=5$. Then there exists a pair of real-valued functions $(w,\xi)\in M$  such that\\
\\[0.3cm]
(1)\quad$S(w,\xi)=m=\inf\limits_{(u,v)\in M}S(u,v)$;\\[0.3cm]
(2)\quad$(w,\xi)$ is a ground state solution to system (2.1).}
\\[0.3cm]
{\bf Remark 2.3.}\quad
 According to the results in Coffman \cite{Co1972}, Kwong \cite{Kw1989}, Chen et al. \cite{CL1991} and McLeod \cite{Mc1993}, the ground state solution of a single nonlinear elliptic equation is unique if the nonlinear term satisfies certain properties. Nevertheless, the uniqueness of ground state solution for system (2.1) remains open.\hfill$\Box$
\\[0.3cm]
{\bf Remark 2.4.}\quad
Recalling formulas (1.3), (1.4) and (2.3), it is easy to check the following relation:\\
$$
\left.
\begin{array}{lll}
&\displaystyle S(\phi,\psi)
+a_2m_2\int_{\mathbb{R^N}}\phi_t^2 dx+\frac{a_1m_1}{2}\int_{\mathbb{R^N}}\psi_t^2 dx
\\[0.3cm]
&\quad=\displaystyle E(\phi ,\psi ,\phi_t,\psi_t)
\\[0.3cm]
&\quad=\displaystyle E(\phi_0,\psi_0,\phi_1,\psi_1)
\\[0.3cm]
&\displaystyle\quad =S(\phi_0,\psi_0)
+a_2m_2\int_{\mathbb{R^N}}\phi_1^2 dx+\frac{a_1m_1}{2}\int_{\mathbb{R^N}}\psi_1^2 dx.
 \end{array}
\right.\eqno(2.10)
$$
\\
\indent In the study of a single Klein-Gordon equation, the ground state solutions of the corresponding elliptic equation play an important role in investigating certain qualitative properties of solutions to its Cauchy problem. In particular, by virtue of the variational characteristics of such ground state solutions, sharp conditions for the finite time blow-up and global existence of solutions to the Cauchy problem have been established (see \cite{Le1974}, \cite{PS1975}). Nevertheless, for the nonlinear Klein-Gordon system (1.1), we are concerned with the following problem.
\\[0.3cm]
{\bf Problem 3:}\quad
Can the ground state solutions of the elliptic system (2.1) be also applied to the research on qualitative properties of solutions to the Cauchy problem associated with system (1.1)? In particular, can we utilize the variational features of ground state solutions for system (2.1) to establish a sharp threshold for between the finite time blow-up and global existence of solutions to the Cauchy problem (1.1)-(1.2)?
\\[0.3cm]
\indent The following theorem gives a positive answer to Problem 3.
\\[0.3cm]
{\bf Theorem 2.3 }.\quad {\sl Let $N=5, m_1>0,m_2>0,a_1>0,a_2>0$.
Suppose $(\phi_0,\psi_0,\phi_1,\psi_1) \in H^1 (\mathbb{R^N})\times H^1(\mathbb{R^N})\times L^{2}(\mathbb{R^N})\times L^{2}(\mathbb{R^N})$ and satisfies
$$
E(\phi_0,\psi_0,\phi_1,\psi_1)<m.\eqno(2.11)
$$
Then the following conclusions hold:
\\[0.3cm]
(1)\quad  If $Q(\phi_0,\psi_0)<0$, then the solution $(\phi,\psi)$ to the Cauchy problem (1.1)-(1.2) blows up in finite time.
\\[0.3cm]
(2)\quad  If $Q(\phi_0,\psi_0)>0,$ then the solution $(\phi,\psi)$ to the Cauchy problem (1.1)-(1.2) exists globally on $[0,\infty)$.}
\\[0.3cm]
{\bf Remark 2.5.}\quad In particular, if $E(\phi_0,\psi_0,\phi_1,\psi_1)<0$,
 then the solution $(\phi,\psi)$ to the Cauchy problem (1.1)-(1.2) blows up in finite time.
 \\[0.3cm]
\indent {\sl Indeed, it follows from $E(\phi_0,\psi_0,\phi_1,\psi_1)<0$ and the expression for $Q(u,v)$ (see (2.4)) that\\
$$ Q(\phi_0,\psi_0)=2E(\phi_0,\psi_0,\phi_1,\psi_1)-2a_1a_2m_1m_2\int_{\mathbb{R^N}}\phi_0\psi_0^2 dx
$$
$$\qquad-2a_2m_2\int_{\mathbb{R^N}}\phi_1^2-a_1m_1\int_{\mathbb{R^N}}\psi_1^2<0,$$
\\[0.2cm]
which yields  the finite time blow-up of solutions by conclusion (1) of Theorem 2.3.}\hfill $\Box$
\\[0.3cm]
\indent Finally, motivated by the studies on the instability of standing wave solutions for the single Klein-Gordon equation (see \cite{BC1981}), we focus on the following problem concerning standing wave solutions with ground states for the nonlinear Klein-Gordon system (1.1).
\\[0.3cm]
 {\bf Problem 4:}\quad Are the standing wave solutions with ground
 states for system (1.1) stable?
 \\[0.3cm]
\indent The following theorem gives a qualitative answer to Problem 4: the standing wave solutions with ground states of system (1.1) are strongly unstable.
\\[0.3cm]
\textbf{Theorem 2.4}.\quad {\sl Let $N=5,~ m_1>0,~m_2>0,~a_1>0,~a_2>0 $, and let $(w,\xi)$ be a ground state solution of the elliptic system (2.1). Then for any $\epsilon>0$, there exists $(\phi_0,\psi_0)\in H^1(\mathbb{R^N})\times H^1(\mathbb{R^N})$ such that\\
$$\parallel\phi_0-w\parallel_{H^1(\mathbb{R^N})}<\epsilon,\quad \parallel\psi_0-\xi\parallel_{H^1(\mathbb{R^N})}<\epsilon.\eqno(2.12)
$$
\\
Moreover, the solution $(\phi,\psi)$ to system (1.1) with the initial conditions
\\[0.3cm]
$$
\left\{
\begin{array}{lll}
 		\phi(0,x)=\phi_0,\; \psi(0,x)=\psi_0,
 \\[0.3cm]
		\phi_t(0,x)=0,\; \psi_t(0,x)=0
\end{array}
\right.\eqno(2.13)
   $$
   \\
is defined on some interval $[0,T)$ with $T\in(0,+\infty)$, satisfying  $(\psi,\phi)\in C([0,T),H^1(\mathbb{R^N})\times H^1(\mathbb{R^N}))$ and
$$
	\lim_{t\to T}\left( \parallel\phi\parallel _{H^1(\mathbb{R^N})}+ \parallel\psi\parallel _{H^1(\mathbb{R^N})}\right)=+\infty.\eqno(2.14)
$$}
\\[0.3cm]
{\bf Remark 2.6.}\quad Recalling Theorem 1.4 in Miyazaki \cite{Miy2021}, for $N=5$ and $\omega=0$, there is a perfect coincidence between our result (Theorem 2.4) and that stated in Theorem 1.4 from \cite{Miy2021}.\hfill$\Box$
\\[0.3cm]
{\bf Remark 2.7.}\quad Theorem 2.4 establishes the instability of standing waves with  ground state for system (1.1). During the proof of this theorem, we can choose the initial data $(\phi_0,\psi_0)$ to be compactly supported. This indicates that such instability of standing wave solutions is a local property rather than a global one.\hfill$\Box$
\section{Finite Time Blow-up}
  \indent In this section, we focus on the finite time blow-up of solutions to the Cauchy problem (1.1)-(1.2). If the initial data $(\phi_0,\psi_0)$ and $(\phi_1,\psi_1)$
 satisfy given conditions and the initial energy $E(\phi_0,\psi_0,\phi_1,\psi_1)$ is negative, then the solution to the Cauchy problem (1.1)-(1.2) blows up in finite time, which is stated in Theorem 2.1. We mainly adopt the concavity analysis proposed in \cite{Le1974} and define a suitable non-negative auxiliary function $F(t)$
 to prove Theorem 2.1.\\
 \indent By Proposition 1.1, for arbitrary initial data $ (\phi_0,\psi_0)\in H^1(\mathbb{R^N})\times H^1(\mathbb{R^N}), ~~(\phi_1,\psi_1)\in L^{2}(\mathbb{R^N})\times L^{2}(\mathbb{R^N})$, the Cauchy problem (1.1)-(1.2) is locally well-posed, and its solutions admit the energy conservation law. On the other hand, a negative initial energy serves as an important guarantee for the finite time blow-up for solutions to the Cauchy problem (1.1)-(1.2).
 \\[0.3cm]
 \indent We now proceed to prove Theorem 2.1.
 \\[0.3cm]
{\bf Proof of Theorem 2.1 }.\quad Let $(\phi,\psi,\phi_{t},\psi_{t})\in C([0,T);\;H^1(\mathbb{R^N})\times H^1(\mathbb{R^N})\times L^{2}(\mathbb{R^N})\times L^{2}(\mathbb{R^N}))$ be the solution of the Cauchy problem (1.1)-(1.2). For any given constant $\alpha\in (0,\frac{1}{2})$, define\\
\bea F(t)=\left(a_2m_2\int_{\mathbb{R^N}}\phi^2dx+\frac{a_1m_1}{2}
\int_{\mathbb{R^N}}\psi^2dx\right)^{-\frac{\alpha}{2}}
:=A^{-\frac{\alpha}{2}}(t),
\eea
\\
where
$$A(t)=a_2m_2\int_{\mathbb{R^N}}\phi^2dx+\frac{a_1m_1}{2}
\int_{\mathbb{R^N}}\psi^2dx.\eqno(3.1a)$$
\\
Then we have
\bea F^{'}(t)=-\frac{\alpha}{2}A^{-\frac{\alpha}{2}-1}
\left(2a_2m_2\int_{\mathbb{R^N}}\phi\phi_tdx+
a_1m_1\int_{\mathbb{R^N}}\psi\psi_tdx\right),
\eea
and
\bea
 F^{''}(t)&\displaystyle=(-\alpha)(-\frac{\alpha+2}{2})A^{-\frac{\alpha+4}{2}}
 2\left(a_2m_2\int_{\mathbb{R^N}}\phi\phi_tdx+\frac{a_1m_1}{2}
 \int_{\mathbb{R^N}}\phi\phi_tdx\right)^2
 \nn\\
  &\displaystyle\qquad -A^ {-\frac{\alpha+2}{2}}\alpha\left(a_2m_2\int_{\mathbb{R^N}}\phi_t^2dx +\frac{a_1m_1}{2}\int_{\mathbb{R^N}} \psi_t^2dx\right.
   \nn\\
  &\displaystyle\qquad\qquad\left.+a_2m_2\int_{\mathbb{R^N}}\phi\phi_{tt}dx+\frac{a_1m_1}{2}
  \int_{\mathbb{R^N}}\psi\psi_{tt}dx\right)
 \nn\\
  	&\displaystyle=A^{-\frac{\alpha+4}{2}}\left[\alpha(\alpha+2)\left(a_2m_2
  \int_{\mathbb{R^N}}\phi\phi_tdx+\frac{a_1m_1}{2}
  \int_{\mathbb{R^N}}\psi\psi_tdx\right)^2\right.
 \nn\\
  &\qquad\displaystyle-A\alpha\left(a_2m_2\int_{\mathbb{R^N}}\phi_t^2dx+\frac{a_1m_1}{2}
  \int_{\mathbb{R^N}}\psi_t^2dx\right)
 \nn\\		
  &\qquad\displaystyle\left.-A\alpha\left(a_2m_2\int_{\mathbb{R^N}}\phi\phi_{tt}dx+\frac{a_1m_1}{2}
  \int_{\mathbb{R^N}}\psi\psi_{tt}dx\right)\right]
  \nn\\
  &\leq \displaystyle A^{-\frac{\alpha+4}{2}}\left\{(\alpha^2+\alpha) \left(a_2m_2\int_{\mathbb{R^N}}
\phi_t^2dx+\frac{a_1m_1}{2}\int_{\mathbb{R^N}}\psi_t^2dx\right)\right.
\nn\\
&\qquad\displaystyle-A\alpha\left[a_2m_2\int_{\mathbb{R^N}}
\phi(\Delta\phi-m_1^2\phi+2a_1m_1\phi\psi)dx\right.
\nn\\
&\qquad\qquad\displaystyle\left.\left.+\frac{a_1m_1}{2}\int_{\mathbb{R^N}}\psi(\Delta\psi-m_2^2\psi
+2a_2m_2\phi^2)dx\right]\right\}.
\nn\\
& \displaystyle\leq A^{-\frac{\alpha+2}{2}}\left[(\alpha^2+\alpha)\left(a_2m_2\int_{\mathbb{R^N}}\phi_t^2dx
+\frac{a_1m_1}{2}\int_{\mathbb{R^N}}\psi_t^2dx\right)\right.
\nn\\
&\qquad\displaystyle+\alpha a_2m_2\int_{\mathbb{R^N}}|\nabla\phi|^2dx+\alpha a_2m_2m_1^2\int_{\mathbb{R^N}}\phi^2dx
\nn\\
&\displaystyle\qquad-2\alpha a_1a_2m_1m_2\int_{\mathbb{R^N}}\phi^2\psi dx+\alpha \frac{a_1m_1}{2}\int_{\mathbb{R^N}}|\nabla\psi|^2dx
\nn\\
& \displaystyle\qquad+\left.\alpha \frac{a_1m_1}{2}m_2^2\int_{\mathbb{R^N}}\psi^2dx-\alpha a_1a_2m_1m_2\int_{\mathbb{R^N}}\phi^2\psi dx\right]
\nn\\
&\displaystyle=A^{-\frac{\alpha+2}{2}}\alpha\left[(\alpha+1)E(\phi,\psi,\phi_{t},
\psi_{t})
-\alpha a_2m_2\int_{\mathbb{R^N}}|\nabla\phi|^2dx\right.
\nn\\
&\qquad\displaystyle-\alpha\frac{a_1m_1}{2}
\int_{\mathbb{R^N}}|\nabla\psi|^2dx-\alpha a_2m_2m_1^2\int_{\mathbb{R^N}}\phi^2dx\nn\\
&\qquad\displaystyle\left.-\alpha\frac{a_1m_1}{2}m_2^2
\int_{\mathbb{R^N}}\psi^2dx+(2\alpha-1)a_1a_2m_1m_2\int_{\mathbb{R^N}}\phi^2\psi dx\right].
\eea
   \\
Here, the first inequality in (3.3) follows from the H\"{o}lder's inequality and the Cauchy-Schwarz inequality, and the last equality from the expression for the energy functional (1.3).  In view of the expressions (1.3) and (1.4),  for $E(\phi_0,\psi_0,\phi_1,\psi_1)=E_{0}<0$,  we have\\
  $$2a_1a_2m_1m_2\int_{\mathbb{R^N}}\phi^2\psi dx>0.\eqno(3.4)$$
 \\
 Since $0<\alpha<\frac{1}{2}$, we have $2\alpha-1<0$, this together with $m_1>0, m_2>0, a_1>0, a_2>0$ and (3.4) yields\\
 $$(2\alpha-1)a_1a_2m_1m_2\int_{\mathbb{R^N}}\phi^2\psi dx<0.\eqno(3.5)$$
 \\
 Note that (1.3), (1.4) and (3.5), the following estimate for (3.3) holds:\\
 $$
\left.
\begin{array}{lll}
		(3.3)&\leq A^{-\frac{\alpha+2}{2}}\alpha(\alpha+1)E(\phi,\psi,\phi_t,\psi_t)
\\[0.3cm]
&=\alpha(\alpha+1)A^{-\frac{\alpha+2}{2}}E(\phi_0,\psi_0,\phi_1,\psi_1)
\\[0.3cm]
&=\alpha(\alpha+1)A^{-\frac{\alpha+2}{2}} E_{0}.
\end{array}
\right.\eqno(3.6)
   $$
 \\
 Therefore,
 combining (3.1), (3.2), (3.3) and (3.6) with the assumptions of Theorem 2.1, we obtain\\
 $$
	F(0)=\left(a_2m_2\int_{\mathbb{R^N}}\phi_0^2dx
+\frac{a_1m_1}{2}\int_{\mathbb{R^N}}\psi_0^2dx\right)^{-\frac{\alpha}{2}}>0,\eqno(3.7)
$$
\\
$$
\left.
\begin{array}{lll}
		F^{'}(0)		&\displaystyle=-\frac{\alpha}{2}A^{-\frac{\alpha}{2}-1}\left(2a_2m_2\int_{\mathbb{R^N}}\phi\phi_tdx
+a_1m_1\int_{\mathbb{R^N}}\psi\psi_tdx\right)|_{t=0}
\\[0.4cm]
		&\displaystyle=-\alpha A^{-\frac{\alpha+2}{2}}(0)\left(a_2m_2\int_{\mathbb{R^N}}\phi_{0}\phi_1dx
+\frac{a_1m_1}{2}\int_{\mathbb{R^N}}\psi_{0}\psi_1dx\right)<0 .
	\end{array}
\right.\eqno(3.8)
   $$
   \\
  Hence, it follows from (3.3), (3.6) and (3.8) that $F(t)$
 is concave with respect to $t$
 provided $E_{0}<0$. Furthermore, combining (3.7) and (3.8), there exists $T^*>0$
 such that $F(t)\rightarrow 0$ as $t\rightarrow T^*$. Together with (3.1), we further deduce $$
	\lim_{t\to T^*} A(t)=\lim_{t\to T^*}\left(a_2m_2\int_{\mathbb{R^N}}\phi^2dx
+\frac{a_1m_1}{2}\int_{\mathbb{R^N}}\psi^2dx\right)=+\infty.
$$
 \indent This completes the proof of Theorem 2.1.\hfill$\Box$\\

 \section{Existence of standing wave with ground state}
 \indent
 Our aim in this section is to prove Theorem 2.2. Here, we shall adopt several methods from the studies on standing wave with ground state for a single nonlinear wave equation \cite{BC1981,St1997,TZ1999}, and combine the variational method, scaling analysis, compactness arguments and the Lagrange multiplier method.\\
\indent Firstly, based on (2.3), (2.4), (2.5) and (2.6), we investigate the variational characterization of ground state solutions for system (2.1).

\subsection{Existence of solutions to the minimization problem (2.6)}
  \quad\\
  \\
\indent In this subsection, we will prove the conclusion (1) of Theorem 2.2.

\indent First of all, we claim the following proposition which presents a useful variational characterization of ground state solutions for system (2.1). In particular, it also provides an alternative approach to prove the existence of ground state solutions stated in Theorem 2.2.\\
 \quad\\
{\bf Proposition 4.1}.  {\sl Assume that the hypotheses in Theorem 2.2 hold. Let
$$M=\{(u,v)\in H^1(\mathbb{R^N})\times H^1(\mathbb{R^N})\setminus \{(0,0)\}\; , Q(u,v)=0\}.\eqno(4.1)$$
 Then the constrained variational problem (2.6) admits a solution. That is, there exists $(w,\xi)\in M$ such that
 $$ S(w,\xi)=\min_{(u,v)\in M}S(u,v).\eqno(4.2)$$}
 \quad\\
 \indent Before proving this proposition, we first present several key priori conclusions.\\
  \quad\\
 {\bf Proposition 4.2}.\quad {\sl Any solution to system (2.1) belongs to the manifold $M$.}
 \\[0.3cm]
 {\bf Proof.}\quad Let a pair of real-valued functions $(u,v)\in  H^1(\mathbb{R^N})\times H^1(\mathbb{R^N})\setminus \{(0,0)\} $ be a solution to system (2.1). Then $(u,v)$ satisfies\\
  $$
\left\{
\begin{array}{lll}
 	-\Delta u+m_1^2u=2a_1m_1uv,&\;(4.3a)
 \\[0.4cm]		
-\Delta v+m_2^2v=2a_2m_2u^2.&\;(4.3b)
\end{array}
\right.\eqno(4.3)
   $$
   \\
 Multiplying both sides of (4.3a) by $2a_2m_2u$, and of (4.3b) by $a_1m_1v$, then integrating by parts with respect to $x$ on $\mathbb{R^N}$ yield \\
   $$
\left.
\begin{array}{lll}
 		\displaystyle 2a_2m_2\int_{\mathbb{R^N}}|\nabla u|^2dx+2a_2m_1^2m_2\int_{\mathbb{R^N}}|u|^2dx-4a_1a_2m_1m_2\int_{\mathbb{R^N}}u^2vdx
 \\[0.4cm]
		 \displaystyle  \quad+a_1m_1\int_{\mathbb{R^N}}|\nabla v|^2dx+a_1m_1m_2^2\int_{\mathbb{R^N}}|v|^2dx-2a_1a_2m_1m_2\int_{\mathbb{R^N}}u^2vdx=0.
\end{array}
\right.
 $$
 \\
  Simplifying this equality leads to $Q(u,v)=0$, which implies $(u,v)\in M$.\hfill$\Box$\\
   \quad\\
 {\bf Proposition 4.3}.\quad {\sl The functional $S(u,v)$ is bounded below on $M$.}\\

   {\bf Proof.}\quad Let $(u,v)\in  H^1(\mathbb{R^N})\times H^1(\mathbb{R^N})\setminus \{(0,0)\} $. It follows from (2.3), (2.4) and (2.5) that the following identity holds on $M$:\\
   $$S(u,v)=\frac{1}{3}\int_{\mathbb{R^N}}\left(a_2m_2|\nabla u|^2+\frac{1}{2}a_1m_1|\nabla v|^2+a_2m_2m_1^2 u ^2+\frac{a_1m_1m_2^{2}}{2} v ^2\right)dx.\eqno(4.3a)
$$
\\
  Notice that $a_1m_1>0, a_2m_2>0$, consequently, $S(u,v)>0$ on the set $M$.\hfill$\Box$ \\
  \quad\\
  {\bf Proposition 4.4}.\quad {\sl Suppose $N=5$ and $a_1>0,a_2>0,m_1>0,m_2>0$. For
$$(u,v)\in  H^1(\mathbb{R^N})\times H^1(\mathbb{R^N})\setminus \{(0,0)\},~~\lambda>0,$$
let
$$
	u_\lambda(x)=\lambda u(x),\quad v_\lambda(x)=\lambda v(x),\eqno(4.4)
$$
\\
then there exists a unique constant $\beta>0$ depending on $ (u,v) $ such that $Q(u_\beta ,v_\beta )=0$. Furthermore, we have\\
\begin{equation}
	\begin{cases}
		Q(u_\lambda,v_\lambda)>0,  &\ \lambda\in(0,\beta);
\\[0.3cm]
		Q(u_\lambda,v_\lambda)<0,  &\ \lambda\in(\beta,\infty);
\\[0.3cm]
		S(u_\beta,v_\beta)\geq S(u_\lambda,v_\lambda), &\ \forall\lambda>0.
	\end{cases}
	\notag
\end{equation}}
\\[0.3cm]
 {\bf Proof.} From (2.3), (2.4) and (4.4), the expressions of $S(u_\lambda,v_\lambda)$ and $Q(u_\lambda,v_\lambda)$  are listed as follows:\\
  $$
\left.
\begin{array}{lll}
			 S(u_\lambda,v_\lambda)
			&\displaystyle=\lambda^2\int_{\mathbb{R^N}}\left(a_2m_2|\nabla u|^2+\frac{a_1m_1}{2}|\nabla v|^2+a_2m_2m_1^2|u|^2 +\frac{a_1m_1m_2^2}{2}|v|^2\right)dx
\\[0.4cm]
&\displaystyle \qquad -2a_1a_2m_1m_2\lambda^3\int_{\mathbb{R^N}}vu^2dx,
\end{array}
\right.\eqno(4.5)
   $$
   	$$
\left.
\begin{array}{lll}
			Q(u_\lambda,v_\lambda)
			&\displaystyle=\lambda^2\int_{\mathbb{R^N}}\left(2a_2m_2|\nabla u|^2+a_1m_1|\nabla v|^2+2a_2m_2m_1^2|u|^2 +a_1m_1m_2^2|v|^2\right)dx
\\[0.4cm]
&\displaystyle \qquad
-6a_1a_2m_1m_2\lambda^3\int_{\mathbb{R^N}}vu^2dx.
		\end{array}
\right.\eqno(4.6)
   $$
   \\
Notice that $a_1m_1>0, a_2m_2>0$, by Proposition 2.3, there exists $\beta=\beta(u,v)>0$ such that
$$
Q(u_\beta,v_\beta)=0.\eqno(4.7)
		$$
In particular,
$$Q(u_\lambda,v_\lambda)>0\ \text{for} ~~0<\lambda<\beta,\qquad Q(u_\lambda,v_\lambda)<0\ \text{for}~~\lambda>\beta.\eqno(4.7a)$$
 On the other hand, direct computation yields
 $$\frac{d}{d\lambda}S(u_\lambda,v_\lambda)=\lambda^{-1}Q(u_\lambda,v_\lambda).\eqno(4.7b)$$
 Combining with (4.7), we conclude that $S(u_\beta,v_\beta)\geq S(u_\lambda,v_\lambda)$ for any $\lambda>0$.\\
\indent This completes the proof of Proposition 4.4.\hfill$\Box$
\\[0.3cm]
{\bf Remark 4.1}.\quad It follows from (4.5) and (4.7) that $S(\lambda u,\lambda v)$  is concave on $(\beta,+\infty)$. Therefore, note that (4.7b), for any $(u,v)\in  H^1(\mathbb{R^N})\times H^1(\mathbb{R^N})\setminus \{(0,0)\}$  with $\lambda^*>\beta$, we have\\
$$\frac{Q(\lambda^* u,\lambda^* v)}{\lambda^*}<\frac{S(u_\beta,v_\beta)-S(u_\lambda^*,v_\lambda^*)}
{\beta-\lambda^*}<0.\eqno(4.7c)$$
(see the following Figure 1)\\
\begin{figure}[H]
\centering
\begin{tikzpicture}
\def\betaVal{1.2}
\def\lambdaVal{2.8}

\pgfmathsetmacro{\Sbeta}{20 - exp(\betaVal)}
\pgfmathsetmacro{\Slambda}{20 - exp(\lambdaVal)}

\begin{axis}[
    axis lines = middle,
    xlabel = {$\lambda$},
    ylabel = {$S(u_{\lambda},v_{\lambda})$},
    xmin = 0.4, xmax = 3.8,
    ymin = -25, ymax = 48,
    samples = 200,
    xtick={}, ytick={},
    xticklabels={}, yticklabels={},
    xtick style = {draw=none},
    ytick style = {draw=none},
    clip = false,
    width = 10cm,
    height = 7cm,
]
\addplot [thick, solid, domain=0.4:3.8] {20 - exp(x)};

\addplot [only marks, mark=*, mark size=3.5pt] coordinates {(\betaVal,\Sbeta)};
\addplot [only marks, mark=*, mark size=3.5pt] coordinates {(\lambdaVal,\Slambda)};

\addplot [thick, dashed, domain=0.4:4.2] {\Sbeta - (\Sbeta - \Slambda)/(\betaVal - \lambdaVal)*(\betaVal - x)};

\addplot [thick, dashdotted, domain=0.8:4.2] {-exp(\lambdaVal)*(x - \lambdaVal) + \Slambda};

\node at (axis cs:0.5,-1) [below] {$O$};
\node at (rel axis cs:1.15,0.24) {$l_1$};
\node at (rel axis cs:1.15,0.07) {$l_2$};
\node at (axis cs:\betaVal,\Sbeta) [below left, xshift=15pt, font=\small] {$(\beta, S(u_{\beta},v_{\beta}))$};
\node at (axis cs:\lambdaVal,\Slambda) [above right, xshift=-15pt,yshift=5, font=\small] {$(\lambda^*, S(u_{\lambda^*},v_{\lambda^*}))$};
\end{axis}
\end{tikzpicture}
\vspace{3mm}

$$\hbox{Figure 1}$$

\label{fig:curve_secant_tangent}
\end{figure}
By Figure 1, (4.7c) is equivalent to $k_{l_{2}}<k_{l_{1}}$. Here, $l_{1}$ denotes the straight line passing through points $(\beta, S(u_\beta,v_\beta))$ and $(\lambda^*, S(u_{\lambda^*},v_{\lambda^*}))$, $l_{2}$ denotes the tangent line passing through point $(\lambda^*, S(u_{\lambda^*},v_{\lambda^*}))$, $k_{l_{1}}$ and $k_{l_{2}}$ represent the slopes of straight lines $l_{1}$ and $l_{2}$, respectively.
Hence, recalling that $(u_\beta,v_\beta)=(\beta u,\beta v)\in M$, if $(u,v)\in H^1(\mathbb{R}^N)\times H^1(\mathbb{R}^N)$ is such that $Q(u,v)<0$,  (4.7), (4.7c) and $\lambda^* >\beta$ give that $\beta\in (0,1)$ and \\
$$Q(u,v)<-\left\{\inf_{(u^*,v^*)\in M}S(u^*,v^*)-S(u,v)\right\}.\eqno \Box$$\\
\quad\\
\indent Now we turn to prove Proposition 4.1 by solving the constrained variational problem (2.6).\\
\indent By Proposition 4.3, Proposition 2.3 and its proof, let $\{(u_n,v_n), n\in\mathbb{N}\}\subset M$ be a minimizing sequence for the constrained variational problem (2.6), that is, $(u_n,v_n)\neq (0,0)$ in the sense of $H^{1}(\mathbb{R^N})\times H^{1}(\mathbb{R^N})$, and as $n\to +\infty$,\\
$$S(u_n,v_n)\to \inf_{(u,v)\in M}S(u,v),\eqno(4.8)$$
\\
$$Q(u_n,v_n)=0.\eqno(4.9)$$
\\
  According to Definition 2.1 and Proposition 2.1, let $u_n^*$ and $v_n^*$  denote the spherical Schwarz symmetrizations of the positive functions $u_n>0$ and $v_n>0$, respectively. Then $u_n^*$ and $v_n^*$ are spherically symmetric and non-increasing functions with respect to the radial variable $r=|x|$.
Moreover, $u_n^*$  and $v_n^*$ satisfy the following properties:\\
$$[(u_{n})_{\lambda}]^*=(u_n^*)_\lambda,\quad [(v_{n })_{\lambda}]^*=(v_n^*)_\lambda,\eqno(4.10)$$
\\
$$\int_{\mathbb{R^N}}|\nabla u_n^*|^2dx\leq \int_{\mathbb{R^N}}|\nabla u_n|^2dx,~~
 \int_{\mathbb{R^N}}|\nabla v_n^*|^2dx\leq \int_{\mathbb{R^N}}|\nabla v_n|^2dx,
\eqno(4.10a)$$
\\
$$ \int_{\mathbb{R^N}}| u_n^*|^\sigma dx= \int_{\mathbb{R^N}}| u_n|^\sigma dx,~~
 \int_{\mathbb{R^N}}|v_n^*|^\sigma dx= \int_{\mathbb{R^N}}| v_n|^\sigma dx, ~\sigma>1.\eqno(4.10b)$$
 \\
Here,
 $$(u_{n})_{\lambda}=\lambda u_{n},~~ (v_{n})_{\lambda}=\lambda v_{n}.$$
 \\
 Now, for the minimizing sequence \(\{(u_n,v_n),\ n\in\mathbb{N}\}\), set
$$w_n=(u_n^*)_{\beta_n},\quad\xi_n=(v_n^*)_{\beta_n},\eqno(4.11)$$
\\
where \(\beta_n>0\) is uniquely determined by the relation\\
$$
	Q(w_n,\xi_n)=Q((u_n^*)_{\beta_n},(v_n^*)_{\beta_n})=0.
	\eqno(4.12)
$$
\\
\indent On the other hand, it follows from (4.10) and (4.11) that\\
$$
	w_n=[(u_n)_{\beta_n}]^*, \; \xi_n=[(v_n)_{\beta_n}]^*.
	$$
\\
In view of Proposition 4.4, together with (4.3a), (4.9),(4.10), (4.10a), (4.10b) and (4.12), we deduce that\\
$$
	 S(w_n,\xi_n) \leq S((u_n)_{\beta_n},(v_n)_{\beta_n}) \leq S(u_n,v_n).
	 \eqno(4.13)
$$
\\
Combining this inequality with (4.12) yields $\{(w_n,\xi_n),\, n\in\mathbb{N}\}\subset M$ and $S(w_n,\xi_n)\leq S(u_n,v_n)$. Consequently, the pair $(w_n,\xi_n)=\left(\left[(u_n)_{\beta_n}\right]^*,\left[(v_n)_{\beta_n}\right]^*\right)$ also forms a minimizing sequence for problem (2.6). Moreover, the sequences $\{w_n\}$ and $\{\xi_n\}$ satisfying $w_n>0$, $\xi_n>0$, are radially symmetric nonincreasing functions of $|x|$.
\\[0.3cm]
\indent Furthermore, we claim:
\\[0.3cm]
\textbf{Lemma 4.1}.\quad  {\sl Suppose $\{(u_n,v_n), n\in\mathbb{N}\}\subset M$  is a sequence such that $S(u_n,v_n)$ is bounded above. Then $(u_n,v_n)$ is bounded in $H^1(\mathbb{R^N})\times H^1(\mathbb{R^N})$.}\\
{\bf Proof.}\quad By (4.3a) and (4.8), there exists $C>0$ such that\\
	$$||u_n||_{H^1(\mathbb{R^N})}^2+||v_n||_{H^1(\mathbb{R^N})}^2 \lesssim S(u_n,v_n)\leq C.\eqno \Box$$\\
With Lemma 4.1, it follows from (4.12) and (4.13) that\\
$$||w_n||_{H^1(\mathbb{R^N})}^2+||\xi_n||_{H^1(\mathbb{R^N})}^2\leq C.\eqno(4.14)$$
\\
Applying the compactness Lemma 2.1, there exists a subsequence $\{w_{n_k}, k\in\mathbb{N}\}\subset\{w_n,n\in\mathbb{N}\}$ such that as $k\to+\infty$,
$$w_{n_k}\rightharpoonup w_\infty \quad \text{weakly~in}~H^1(\mathbb{R^N}).\eqno(4.14a)$$
\\
On the other hand, for $\{\xi_{n_k}, k\in\mathbb{N}\}\subset\{\xi_n,n\in\mathbb{N}\}$, there exists a subsequence $\{\xi_{n_{k_m}}, m\in\mathbb{N}\}\subset\{\xi_{n_k}, k\in\mathbb{N}\}$ such that\\
$$
\left\{
\begin{array}{lll}
	\xi_{n_{k_m}}\rightharpoonup\xi_\infty \quad \text{weakly~in}~	 H^1(\mathbb{R^N}),\\[0.3cm]
		\xi_{n_{k_m}}\rightarrow\xi_\infty \quad a.e.\ \text{in}\ \mathbb{R^N}.
\end{array}
\right.\eqno(4.14b)
   $$
   \\
Naturally, it follows from (4-14a) that there exists a subsequence $\{w_{n_{k_m}}, m\in\mathbb{N}\}\subset\{w_{n_k}, k\in\mathbb{N}\}$ such that\\
 $$
\left\{
\begin{array}{lll}
	w_{n_{k_m}}\rightharpoonup w_\infty \quad \text{weakly~in}~	 H^1(\mathbb{R^N}),\\[0.3cm]
		w_{n_{k_m}}\rightarrow w_\infty \quad a.e.\  \text{in}\ \mathbb{R^N}.
\end{array}
\right.\eqno(4.14c)
   $$
   \\
   Without ambiguity, we still use $\{(w_n,\xi_n),n\in\mathbb{N}\}$ to denote $\{(w_{n_{k_m}},\xi_{n_{k_m}}), m\in\mathbb{N}\}$.\\
\indent Recalling Lemma 2.1 that the embedding $H_r^1(\mathbb{R}^N)\hookrightarrow L^\sigma(\mathbb{R}^N)$ is compact provided $2<\sigma<\dfrac{2N}{N-2}$, combining (4-14b) and (4-14c) together leads to as $n\to\infty$,\\
$$
\left\{
\begin{array}{lll}
	w_n\rightharpoonup w_\infty \quad \text{strongly~in}~	 L^{2q}(\mathbb{R^N}),\\[0.3cm]
		\xi_n\rightharpoonup \xi_\infty\quad \text{strongly~in}~	 L^{p}(\mathbb{R^N}),
\end{array}
\right.\eqno(4.15)
   $$
   \\
 where $p$ and $q$ satisfy\\
$$\frac{1}{p}+\frac{1}{q}=1,\quad N=5,\quad
        2+\frac{4}{N}<p<\frac{2N}{N-2},\quad
		\frac{2N}{N+2}<q<\frac{2N+4}{N+4}<\frac{N}{N-2}.
\eqno(4.15a)
   $$
   \\
   Accordingly, we assert that\\
   $$(w_\infty,\xi_\infty)\neq(0,0).\eqno(4.16)$$
   \\
   \indent We prove (4.16) by contradiction in what follows.
   \\[0.3cm]
   \indent Assume
   $$(w_\infty,\xi_\infty)=(0,0).\eqno(4.16a)$$
   \\
Since $(w_n,\xi_n)\in M$, it follows from $Q(w_n,\xi_n)=0$ that\\
  $$
\left.
\begin{array}{lll}
		&\displaystyle 2a_2m_2\int_{\mathbb{R^N}}|\nabla w_n|^2dx+a_1m_1\int_{\mathbb{R^N}}|\nabla\xi_n|^2dx
 \\[0.3cm]		
 &\qquad\displaystyle +2a_2m_1^2m_2\int_{\mathbb{R^N}}w_n^2dx+a_1m_2^2m_1\int_{\mathbb{R^N}}
 \xi_n^2dx\\[0.3cm]
		&\quad\displaystyle =6a_1a_2m_1m_2\int_{\mathbb{R^N}}w_n^2\xi_ndx.\end{array}
\right.\eqno(4.16b)
   $$
   \\
  On the other hand, for $\dfrac1p+\dfrac1q=1$, the Young's inequality gives that\\
  $$\int_{\mathbb{R^N}}w_n^2\xi_ndx\leq
   \int_{\mathbb{R^N}}\frac{|\xi_n|^p}{p}dx+\int_{\mathbb{R^N}}\frac{|w_n|^{2q}}{q}
   dx
		\leq C\left(||\xi_n||_{L^p}^p+||w_n||_{L^{2q}}^{2q}\right).\eqno(4.16c) $$
\\
Further, observing (4.15),(4.15a) and (4.16a), as $n\to\infty$,\\
$$
\left\{
\begin{array}{lll}
	w_n\rightharpoonup 0\quad \text{strongly~in}~	 L^{2q}(\mathbb{R^N}),\\[0.3cm]
		\xi_n\rightharpoonup 0 \quad \text{strongly~in}~	 L^{p}(\mathbb{R^N}).
\end{array}
\right.\eqno(4.16d)
   $$
   \\
   Combining this with (4.16b) and (4.16c), it follows that as $n\to\infty$,\\

   $$
\left.
\begin{array}{lll}
		\displaystyle &2a_2m_2\int_{\mathbb{R^N}}|\nabla w_n|^2dx+a_1m_1\int_{\mathbb{R^N}}|\nabla\xi_n|^2dx
\\[0.3cm]
		&\displaystyle\qquad+2a_2m_1^2m_2\int_{\mathbb{R^N}} u^2dx+a_1m_2^2m_1\int_{\mathbb{R^N}}v^2dx-6a_1a_2m_1m_2\int_{\mathbb{R^N}}vu^2dx.
\end{array}
\right.\eqno(4.16e)
   $$
   \\
  With this, we can deduce that as $n\to\infty$,\\
$$2a_2m_2\int_{\mathbb{R^N}}|\nabla w_n|^2dx+a_1m_1\int_{\mathbb{R^N}}|\nabla\xi_n|^2dx\to0.\eqno(4.16f)$$
\\
Moreover, from (4.16b), (4.16c), the Gagliardo-Nirenberg inequality and the Cauchy-Schwarz inequality, we conclude that\\
$$
\left.
\begin{array}{lll}
		&\displaystyle 2a_2m_2\int_{\mathbb{R^N}}|\nabla w_n|^2dx+a_1m_1\int_{\mathbb{R^N}}|\nabla\xi_n|^2dx\\[0.4cm]
		&\displaystyle\qquad \leq6a_1a_2m_1m_2\int_{\mathbb{R^N}}w_n^2\xi_ndx\\[0.4cm]
		&\displaystyle \qquad \leq C\left(||\xi_n||_{L^p(\mathbb{R^N})}^p+||w_n||_{L^{2q}(\mathbb{R^N})}^{2q}\right)\\[0.4cm]
		&\displaystyle\qquad  \leq C||\xi_n||_{L^2(\mathbb{R^N})}^{p-\frac{N}{2}(p-2)}||\nabla\xi_n||_{L^2(\mathbb{R^N})}
^{\frac{N}{2}(p-2)}+C||w_n||_{L^2(\mathbb{R^N})}^{2q-\frac{N}{2}(2q-2)}||\nabla w_n||_{L^2(\mathbb{R^N})}^{\frac{N}{2}(2q-2)}.
  \end{array}
\right.\eqno(4.17)
$$
\\
Direct computaion from (4.15a) gives\\
$$
\left\{
\begin{array}{lll}
		\displaystyle\frac{N}{4}(p-2)>1, \frac{4N}{N+2}<2q<\frac{4N+8}{N+4},
\\[0.5cm]		
\displaystyle\frac{N}{4}(2q-2)>\frac{N}{4}\left(\frac{4N}{N+2}-2\right)
=\frac{N}{4} \frac{2N-4}{N+2} =\frac{N}{2} \frac{N-2}{N+2} >1.
\end{array}
\right.\eqno(4.17a)
   $$
 \\[0.3cm]
Let $\displaystyle\theta=\max\left\{\frac{N}{4}(p-2), \frac{N}{4}(2q-2)\right\}>1$. then (4.14) implies that\\
$$||w_n||_{L^2(\mathbb{R^N})}\leq C, \quad ||\xi_n||_{L^2(\mathbb{R^N})}\leq C.\eqno(4.17b)$$
  \\
 Therefore, it can be inferred from (4.17), (4.17a) and (4.17b) that \\
 $$
\left.
\begin{array}{lll}
		&\displaystyle 2a_2m_2\int_{\mathbb{R^N}}|\nabla w_n|^2dx+a_1m_1\int_{\mathbb{R^N}}|\nabla\xi_n|^2dx
\\[0.4cm]
		&\displaystyle \qquad\leq C\left[2a_2m_2\int_{\mathbb{R^N}}|\nabla w_n|^2dx+a_1m_1\int_{\mathbb{R^N}}|\nabla\xi_n|^2dx\right]^\theta.
\end{array}
\right.
   $$
   \\
   This implies for $\theta>1$ that\\
$$2a_2m_2\int_{\mathbb{R^N}}|\nabla w_n|^2dx+a_1m_1\int_{\mathbb{R^N}}|\nabla\xi_n|^2dx\geq C>0.\eqno(4.17c)$$\\
This contradicts (4.16f). Hence, (4.16) holds, i.e., $(w_\infty,\xi_\infty)\neq(0,0)$.

Next, set
$$w=(w_\infty)_\mu, \quad \xi=(\xi_\infty)_\mu,\eqno(4.18)$$
\\
where $\mu>0$ is determined by
$$Q(w,\xi)=Q((w_\infty)_\mu,(\xi_\infty)_\mu).\eqno(4.18a)$$
\\
Combining (4.14b), (4.14c), (4.15) and (4.18), we conclude that as $n\to\infty$, \\
$$
\left\{
\begin{array}{lll}
	(w_n)_\mu\rightarrow w  ~~\text{strongly~~in}~L^{2q}(\mathbb{R^N}),
\\[0.3cm]
(\xi_n)_\mu\rightarrow \xi ~~\text{strongly~~in}~L^{2q}(\mathbb{R^N}),
\\[0.3cm]
		(w_n)_\mu\rightharpoonup w  \quad (\xi_n)_\mu\rightharpoonup \xi~~\text{weakly~~in}~H^{1}(\mathbb{R^N}),
\end{array}
\right.\eqno(4.19)
   $$
 \\
 where $p$ and $q$ are determined by (4.15a). Recalling that $Q(w_n,\xi_n)=0$,
  we deduce from Proposition 4.4 that\\
  $$
S\left[(w_n)_\mu,(\xi_n)_\mu\right]\leq S(w_n,\xi_n).\eqno(4.20)
   $$
   \\
   Thus, using (4.19) and (4.20), we obtain\\
   $$S(w,\xi)\leq\varliminf_{n\to\infty}S[(w_n)_\mu,(\xi_n)_\mu]\leq \lim_{n\to\infty}S(w_n,\xi_n)=\inf_{(u,v)\in M}S(u,v).\eqno(4.21)$$
   \\
   On the other hand, since $(w,\xi)\neq(0,0)$ and $ Q(w,\xi)=0$, we have $(w,\xi)\in M$. Moreover, it follows from (4.21) that $(w,\xi)$
 satisfies
 $$S(w,\xi)=\min_{(u,v)\in M}S(u,v).$$
 Thus Proposition 4.1 is proved. Accordingly, conclusion (1) in Theorem 2.2 holds.
\hfill$\Box$
\subsection{Existence of Ground State for (2.1)}
\quad\\
\\
\indent In this subsection, we prove the conclusion (2) in Theorem 2.2.\\
\indent Since $(w,\xi)$ is a solution to the minimization problem (4.2), there exists a Lagrange multiplier $\Lambda$  such that\\
$$\delta_w[S(w,\xi)+\Lambda Q(w,\xi)]=0, \quad \delta_\xi[S(w,\xi)+\Lambda Q(w,\xi)]=0,\eqno(4.22)$$\\
where $\delta_u G$ denotes the variation of
$G(u,v)$ with respect to $u$.\\
\indent Observe that\\
$$\delta_uG(u,v)=\frac{\partial}{\partial\epsilon}G(u+\epsilon\delta u,v)|_{\epsilon=0},\eqno(4.22a)$$
it follows that\\
$$
\left\{
\begin{array}{lll}
		\displaystyle\delta_w\left[S(w,\xi)+\Lambda Q(w,\xi)\right]\\[0.4cm]
\displaystyle\quad=\left\langle  2(1+2\Lambda)\left(-a_2m_2\Delta w+a_2m_2m_1^2w\right)-4a_1a_2m_1m_2(1+3\Lambda)w\xi,\delta w \right\rangle,
\\[0.5cm]
		\displaystyle\delta_\xi\left[S(w,\xi)+\Lambda Q(w,\xi)\right]\\[0.4cm]
\displaystyle\quad=\left\langle (1+2\Lambda)\left(-a_1m_1\Delta \xi+a_1m_1m_2^2\xi\right)-2a_1a_2m_1m_2(1+3\Lambda)w^2,\delta \xi \right\rangle,
\end{array}
\right. \eqno(4.22b)
   $$
   \\
where $\delta u$ stands for the variation of $u$, and $\displaystyle \left\langle f,~g\right\rangle =\int_{\mathbb{R^N}}fgdx$. Hence, (4.22) is equivalent to the following two equations:\\
$$
\left\{
\begin{array}{lll}
		\displaystyle\int_{\mathbb{R^N}}2(1+2\Lambda)\left(a_2m_2|\nabla w|^2+a_2m_2m_1^2w^2\right)dx
=\int_{\mathbb{R^N}}4a_1a_2m_1m_2(1+3\Lambda)w^2\xi dx,\\[0.5cm]
		\displaystyle\int_{\mathbb{R^N}} (1+2\Lambda)\left(a_1m_1|\nabla\xi|^2+a_1m_1m_2^2\xi^2\right)dx=\int_{\mathbb{R^N}}2a_1a_2m_1m_2(1+3\Lambda)w^2\xi dx.
\end{array}
\right.\eqno(4.22c)
   $$
   \\
Note that $Q(w,\xi)=0$ is of the form\\
  $$
\left.
\begin{array}{lll}
		&\displaystyle \int_{\mathbb{R^N}}(2a_2m_2|\nabla w|^2+a_1m_1|\nabla\xi|^2+2a_2m_2m_1^2w^2+a_1m_1m_2^2\xi^2)dx\\[0.3cm]
		&\displaystyle\qquad=6a_1a_2m_1m_2\int_{\mathbb{R^N}}w^2\xi dx,
\end{array}
\right.\eqno(4.22d)
   $$
\\
this together with (4.22c) yields\\
 $$
\left.
\begin{array}{lll}
		&\displaystyle (1+2\Lambda)\int_{\mathbb{R^N}}6a_1a_2m_1m_2w^2\xi dx=(1+3\Lambda)\int_{\mathbb{R^N}}6a_1a_2m_1m_2w^2\xi dx\\[0.5cm]
		&\displaystyle \qquad=(1+2\Lambda)\int_{\mathbb{R^N}}\left(2a_2m_2|\nabla w|^2+a_1m_1|\nabla\xi|^2 +2a_2m_2m_1^2w^2+a_1m_1m_2^2\xi^2\right)dx.
\end{array}
\right.\eqno(4.22e)
   $$
\\
Note that  $(w,\xi)\neq(0,0)$, there holds ~~ $\Lambda=0$. \\
\indent Thus, combining (4.22), (4.22b) and (4.22e) together leads to\\
$$
\left\{
\begin{array}{lll}
		-a_2m_2\Delta w+a_2m_2m_1^2w=2a_1a_2m_1m_2w\xi,\\[0.4cm]
		-a_1m_1\Delta\xi+a_1m_1m_2^2\xi=2a_1a_2m_1m_2w^2.
\end{array}
\right.\eqno(4.22f)
   $$
   \\
In view of $a_1\neq0, a_2\neq0, m_1\neq0, m_2\neq0$, (4.22f) is equivalent to\\
$$
\left\{
\begin{array}{lll}
		-\Delta w+m_1^2w=2a_1m_1w\xi,\\[0.3cm]
		-\Delta\xi+m_2^2\xi=2a_2m_2w^2.
\end{array}
\right.\eqno(4.22g)
   $$
   \\
 This shows that $(w,\xi)$ is a solution to system (2.1). Recalling (4.2), $(w,\xi)$ is also a ground state solution of system (2.1). Consequently, conclusion (2) in Theorem 2.2 is proved.\\
 \indent Thus, the proof of Theorem 2.2 is completed.\hfill$\Box$\\
 \\
  \textbf{Remark 4.2}.~~ According to Theorem 2.2, $(w,\xi)$ is a ground state solution of (2.1). Then it follows from (2.2) that $(\phi(t,x),\psi(t,x))=(w(x),\xi(x))$ is a standing wave solution to (1.1). Therefore, we establish the existence of standing wave with ground state for system (1.1).\hfill$\Box$
 \section{ Sharp Criterion for Blow-up and Global Existence}
\indent In this section, the potential well analysis proposed in the study of hyperbolic partial differential equations \cite{PS1975} and the concavity analysis \cite{Le1974} are adopted to prove Theorem 2.3. Specifically, the sharp conditions for blow-up and global existence of solutions to the Cauchy problem (1.1)-(1.2) are derived by constructing invariant sets under the flow generated by the problem.\\
\indent By (1.3) and (2.3), the energy functional $E(\phi,\psi,\phi_t,\psi_t)$ can be rewritten as\\
$$
\left\{
\begin{array}{lll}
		&\displaystyle E(\phi,\psi,\phi_{t},\psi_{t})=S(\phi,\psi)+a_2m_2\int_{\mathbb{R^N}}|\phi_{t}|^2dx
+\frac{a_1m_1}{2}
\int_{\mathbb{R^N}}|\psi_{t}|^2dx,
\\[0.5cm]
		&\displaystyle \qquad(\phi,\psi)\in H^1(\mathbb{R^N})\times H^1(\mathbb{R^N}),~~(\phi_{t},\psi_{t})\in L^2(\mathbb{R^N})\times L^2(\mathbb{R^N}).
\end{array}
\right.\eqno(5.1^{*})
   $$
\\
\indent Firstly, we construct two invariant manifolds for the Cauchy problem (1.1)-(1.2).
\\[0.3cm]
{\bf Lemma 5.1}.\quad {\sl For $m$ in (2.6), define\\
$$
\left.
\begin{array}{lll}
		K_1:=\left \{(\phi,\psi)\in H^1(\mathbb{R^N})\times H^1(\mathbb{R^N}):~~ S(\phi,\psi)<m, ~~Q(\phi,\psi)<0\right\},
\end{array}
\right.\eqno(5.1)
   $$
  $$
\left.
\begin{array}{lll}
		K_2:=\left \{(\phi,\psi)\in H^1(\mathbb{R^N})\times H^1(\mathbb{R^N}) :~~S(\phi,\psi)<m, ~~Q(\phi,\psi)>0 \right\}\cup\{(0,0)\}.
\end{array}
\right.\qquad\eqno(5.2)
   $$
   \\
 Then $K_1$ and $K_2$ are invariant under the flow associated with the Klein-Gordon system (1.1) provided $E(\phi_0,\psi_0,\phi_1,\psi_1)<m$. That is, for any $(\phi_0,\psi_0)\in K_i, ~i=1,2$, if $E(\phi_0,\psi_0,\phi_1,\psi_1)<m$, then the solution $(\phi(t),\psi(t))$ of the Cauchy problem (1.1)-(1.2) satisfies $(\phi,\psi)\in K_i, ~i=1,2$ for $t\in [0,T)$, where $T$ denotes the maximal existence time.}
 \\[0.3cm]
{\bf Proof}.\quad
    Let $(\phi_0,\psi_0)\in K_1$,  and let $(\phi(t),\psi(t))$ be the solution of the Cauchy problem (1.1)-(1.2) on $[0,T)$.  From (5.1$^{*}$) and Lemma 1.1 it follows that\\
$$S(\phi,\psi)\leq E(\phi ,\psi ,\phi_t,\psi_t)=E(\phi_0,\psi_0,\phi_1,\psi_1)<m.
\eqno(5.3)$$
\\
In the following, we will check
$$ Q(\phi(t),\psi(t))<0,\quad \forall\ t\in [0,T).\eqno(5.3a)$$
 We will show this by contradiction. Assume (5.3a) is not true and there is a $t^*$ such that $(\phi(t^{*}),\psi(t^{*})) \notin K_{1}$. By lower semi-continuity of $Q(\phi(t),\psi(t))$, there is a minimal $t_{1}$ such that $(\phi(t_{1}),\psi(t_{1})) \notin K_{1}$, i.e.,
 \\
  $$Q(\phi(t_{1}),\psi(t_{1}))=0,\eqno(5.3b)$$
 and
 $$Q(\phi(t),\psi(t))<0\quad \text{for}~~0\leq t<t_{1}.\eqno(5.3c)$$
 \\
 \indent We further claim that
 $$(\phi(t_{1}),\psi(t_{1}))\in H^1(\mathbb{R^N})\times H^1(\mathbb{R^N})\setminus \{(0,0)\}.\eqno(5.3d)$$
 \\
 \indent {In fact, by Lemma 2.3 and its proof, we have for $t>t_{1}$, $0<\|\phi(t)\|_{H^1(\mathbb{R^N})}+\|\psi(t)\|_{H^1(\mathbb{R^N})}<\varepsilon$ with $\varepsilon$ small sufficiently, there holds \\
 $$Q(\phi(t),\psi(t))>0.\eqno(5.3e)$$
 \\
 On the other hand, Sobolev embedding and (5.3c) imply that for $0\leq t<t_{1}$, there exists $\delta>0$ large enough such that\\ $$\|\phi(t)\|_{H^1(\mathbb{R^N})}+\|\psi(t)\|_{H^1(\mathbb{R^N})}\geq \delta>0.\eqno(5.3f)$$
 \\
 Hence, by continuity of $Q(\phi(t),\psi(t))$ and of $\|\phi(t)\|_{H^1(\mathbb{R^N})}+\|\psi(t)\|_{H^1(\mathbb{R^N})}$, (5.3b) yields that  $(\phi(t_{1}),\psi(t_{1}))\in H^1(\mathbb{R^N})\times H^1(\mathbb{R^N})\setminus \{(0,0)\}$. This is (5.3d).}\\

 \indent On the other hand, combining (5.3b) and (5.3d) together derives that
 $(\phi(t_{1}),\psi(t_{1}))\in M$.  It then follows from (2.5) and (2.6) that $S(\phi(t_{1}),\psi(t_{1}))\geq m$, and this contradicts (5.3). Therefore $K_{1}$ is an invariant set. Similarly we can check that $K_{2}$ is also an invariant manifold.\hfill$\Box$\\
\\
\indent Next, we establish a crucial estimate for the  functional $S(u,v)$.\\
\\
  {\bf Lemma 5.2 }.\quad {\sl For $(u,v)\in  H^1(\mathbb{R^N})\times H^1(\mathbb{R^N})\setminus\{(0,0)\},$ ~and $\lambda>0$, let\\
 $$u_\lambda(x)=\lambda u(x), \quad v_\lambda(x)=\lambda v(x).\eqno(5.4)$$\\
 If there exists $0<\alpha<1$ such that
 $$Q(u_\alpha,v_\alpha)=0,\eqno(5.4a)$$
 then there holds\\
 $$S(u,v)-S(u_\alpha,v_\alpha)\geq\frac{1}{2}Q(u,v).\eqno(5.4b)$$}
 \\
{\bf Proof.}\quad In view of (2.3), (2.4) and (5.4), direct computation  yields that\\
$$Q(u_\lambda,v_\lambda)=A\lambda^2-B\lambda^3,\quad S(u_\lambda,v_\lambda)=\frac{1}{2}A\lambda^2-\frac{1}{3}B\lambda^3,\eqno(5.4c)$$
\\
where\\
$$
\left\{
\begin{array}{lll}
			\displaystyle A=2a_2m_2\int_{\mathbb{R^N}}|\nabla u|^2dx+a_1m_1\int_{\mathbb{R^N}}|\nabla v|^2dx
\\[0.4cm]
\displaystyle\qquad+2a_2m_1^2m_2\int_{\mathbb{R^N}}u^2dx+a_1m_2^2m_1
\int_{\mathbb{R^N}}v^2dx,\\[0.4cm]
			\displaystyle B=6a_1a_2m_1m_2\int_{\mathbb{R^N}}u^2vdx.
		\end{array}
\right.\eqno(5.4d)
   $$
   \\
On the other hand, $Q(u_\alpha,v_\alpha)=0$ implies that\\
$$A\alpha^2=B\alpha^3,\eqno(5.4e)$$
which together with (5.4c) yields
$$S(u_\alpha,v_\alpha)=\frac{1}{6}B\alpha^3.\eqno(5.4f)$$
Note that
$$Q(u,v)=A-B, \quad S(u,v)=\frac{1}{2}A-\frac{1}{3}B,$$
\\
it follows from $0<\alpha<1$ that
\bea
S(u,v)-S(u_\alpha,v_\alpha)	=\frac{1}{2}A-\frac{1}{3}B-\frac{1}{6}B\alpha^3>\frac{1}{2}A-\frac{1}{3}B
-\frac{1}{6}B=\frac{1}{2}A-\frac{1}{2}B=\frac{1}{2}Q(u,v).
\eea
  \indent This finishes the proof of Lemma 5.2.\hfill$\Box$\\
  \\
  \indent Now we are in position to prove Theorem 2.3.
  \\[0.3cm]
    {\bf Proof of Theorem 2.3}.  Firstly, we prove conclusion (1) of Theorem 2.3.
    \\[0.3cm]
 {\bf Proof of (1) }.\quad  Let $(\phi,\psi)$ be the solution of the Cauchy problem (1.1)-(1.2) on $[0,T)$. Let\\
  $$F(t)=\int_{\mathbb{R^N}}\left(a_2m_2\phi^2+\frac{a_1m_1}{2}\psi^2\right)dx.
  \eqno(5.5)$$\\
  Direct calculation gives\\
   $$F^{'}(t)=2\int_{\mathbb{R^N}}\left(a_2m_2\phi\phi_t
+\frac{a_1m_1}{2}\psi\psi_t\right)dx
  \eqno(5.5a)$$
and
  $$
\left.
\begin{array}{lll}
		  F^{''}(t)			&\displaystyle=2a_2m_2\int_{\mathbb{R^N}}\phi_t^2dx+2a_2m_2
\int_{\mathbb{R^N}}\phi\phi_{tt}dx\\[0.4cm]
&\displaystyle\qquad+a_1m_1\int_{\mathbb{R^N}}\psi_t^2dx+a_1m_1
\int_{\mathbb{R^N}}\psi\psi_{tt}dx\\[0.4cm]
			&\displaystyle=2a_2m_2\int_{\mathbb{R^N}}\phi_t^2dx
+2a_2m_2\int_{\mathbb{R^N}}\phi(\Delta\phi-m_1^2\phi+2a_1m_1\phi\psi)dx\\[0.4cm]			&\displaystyle\qquad+a_1m_1\int_{\mathbb{R^N}}\psi_t^2dx+a_1m_1
\int_{\mathbb{R^N}}\psi(\Delta\psi-m_2^2\psi+2a_2m_2\phi^2)dx\\[0.4cm]			&\displaystyle=2a_2m_2\int_{\mathbb{R^N}}\phi_t^2dx+2a_2m_2\int_{\mathbb{R^N}}
(-|\nabla\phi|^2-m_1^2\phi^2+2a_1m_1\psi\phi^2)dx\\[0.4cm]			&\displaystyle\qquad+a_1m_1\int_{\mathbb{R^N}}\psi_t^2dx+a_1m_1
\int_{\mathbb{R^N}}(-|\nabla\psi|^2-m_2^2\psi^2+2a_2m_2\psi\phi^2)dx\\[0.4cm]			&\displaystyle=2a_2m_2\int_{\mathbb{R^N}}\phi_t^2dx+a_1m_1
\int_{\mathbb{R^N}}\psi_t^2dx
-2a_2m_2\int_{\mathbb{R^N}}|\nabla\phi|^2dx\\[0.4cm]
&\displaystyle\qquad-a_1m_1\int_{\mathbb{R^N}}|\nabla\psi|^2dx
-2a_2m_2m_1^2\int_{\mathbb{R^N}}\phi^2dx\\[0.4cm]
&\displaystyle\qquad-a_1m_1m_2^2\int_{\mathbb{R^N}}\psi^2dx
+6a_1a_2m_1m_2\int_{\mathbb{R^N}}\psi\phi^2dx\\[0.4cm]		&\displaystyle=2a_2m_2\int_{\mathbb{R^N}}\phi_t^2dx+a_1m_1
\int_{\mathbb{R^N}}\psi_t^2dx-Q(\phi,\psi).
	\end{array}
\right.\eqno(5.5b)
   $$
  \\
  It follows from Lemma 5.1 and $Q(\phi_0,\psi_0)<0$ that $Q(\phi ,\psi )<0$. This together with (5.5b) gives\\
  $$F^{''}(t)>0, \quad t\in [0,T).\eqno(5.5c)$$
  \\
  On the other hand, combining (1.3), (1.4) with (5.5b), we deduce that\\
$$
\left.
\begin{array}{lll}
		F^{''}(t)&\displaystyle=5a_2m_2\int_{\mathbb{R^N}}\phi_t^2dx
+\frac{5}{2}a_1m_1\int_{\mathbb{R^N}}\psi_t^2dx
\\[0.4cm]			&\displaystyle\quad+a_2m_2\int_{\mathbb{R^N}}|\nabla\phi|^2dx+\frac{1}{2}a_1m_1
\int_{\mathbb{R^N}}|\nabla\psi|^2dx\\[0.4cm]
&\displaystyle\quad+a_2m_2m_1^2\int_{\mathbb{R^N}}\phi^2dx
+\frac{a_1m_1}{2}m_2^2\int_{\mathbb{R^N}}\psi^2dx\\[0.4cm]
&\quad-3E(\phi_0,\psi_0,\phi_1,\psi_1).
\end{array}
\right.\eqno(5.5d)
   $$
\\
Furthermore, we claim\\
\\
{\bf Proposition 5.1}.\quad {\sl There exists $t_1\in[0,T)$ such that\\
$$F^{'}(t_1)>0.\eqno(5.5e)$$}
\\
\indent The proof of Proposition 5.1 will be given in the final part of the proof of conclusion (1) of Theorem 2.3.\hfill$\sharp$
\\[0.3cm]
\indent
Note that  F(t)
 is convex by (5.5c),  it follows from (5.5e) that
$F(t)$ is strictly increasing in $t$ for all $t>t_1$
 within the maximal existence interval. Let $m^*=\min\left\{m_1^2, m_2^2\right\}>0$,
 then the quantity\\
 $$
\left.
\begin{array}{lll}		
&\displaystyle a_2m_2m_1^2\int_{\mathbb{R^N}}\phi^2dx+\frac{a_1m_1}{2}m_2^2\int_{\mathbb{R^N}}
\psi^2dx-3E(\phi_0,\psi_0,\phi_1,\psi_1)\\[0.4cm]
			&\displaystyle\geq m^*\left(a_2m_2\int_{\mathbb{R^N}}\phi^2dx
+\frac{a_1m_1}{2}\int_{\mathbb{R^N}}\psi^2dx\right)
-3E(\phi_0,\psi_0,\phi_1,\psi_1)\\[0.4cm]
			&\displaystyle=m^*F(t)-3E(\phi_0,\psi_0,\phi_1,\psi_1)
\end{array}
\right.\eqno(5.5f)
   $$
\\
will eventually become positive, and will remain positive thereafter. Thus for $t$ large enough, from (5.5d) it follows that \\
$$F^{''}(t)\geq5\left(a_2m_2\int_{\mathbb{R^N}}\phi_t^2dx
+\frac{a_1m_1}{2}\int_{\mathbb{R^N}}\psi_t^2dx\right).\eqno(5.5g) $$
\\
Combining (5.5), (5.5a) with (5.5g) yields\\
$$
\left.
\begin{array}{lll}
		&\displaystyle F(t)F^{''}(t)\geq\left(\int_{\mathbb{R^N}}a_2m_2\phi^2dx
+\int_{\mathbb{R^N}}\frac{a_1m_1}{2}\psi^2dx\right)\\[0.5cm]			&\displaystyle\qquad\qquad\qquad\quad\cdot 5\left(a_2m_2\int_{\mathbb{R^N}}\phi_t^2dx
+\frac{a_1m_1}{2}\int_{\mathbb{R^N}}\psi_t^2dx\right).
\end{array}
\right.\eqno(5.5h)
   $$
\\
Furthermore, by virtue of the H\"{o}lder inequality and the Cauchy-Schwarz inequality, we obtain\\
$$
\left.
\begin{array}{lll}
		[F^{'}(t)]^2
			&\displaystyle=\left[2a_2m_2\int_{\mathbb{R^N}}\phi\phi_tdx
+a_1m_1\int_{\mathbb{R^N}}\psi\psi_tdx\right]^2\\[0.4cm]		&\displaystyle\leq\left[2a_2m_2\left(\int_{\mathbb{R^N}}\phi^2dx\right)
^{\frac{1}{2}}
\left(\int_{\mathbb{R^N}}\phi_t^2dx\right)^{\frac{1}{2}}
+a_1m_1\left(\int_{\mathbb{R^N}}\psi^2dx\right)^{\frac{1}{2}}
\left(\int_{\mathbb{R^N}}\psi_t^2dx\right)^{\frac{1}{2}}\right]^2\\[0.4cm]
			&\displaystyle\leq 4\left(a_2m_2\int_{\mathbb{R^N}}\phi^2dx
+\frac{a_1m_1}{2}\int_{\mathbb{R^N}}\psi^2dx\right)\\[0.4cm]
&\displaystyle\qquad\cdot
\left(a_2m_2\int_{\mathbb{R^N}}\phi_t^2dx+\frac{a_1m_1}{2}
\int_{\mathbb{R^N}}\psi_t^2dx\right)\\[0.4cm]
			&\displaystyle\leq\frac{4}{5}F(t)F^{''}(t),
\end{array}
\right. \\[0.3cm]
   $$
   \\
   which is equivalent to \\
   $$F(t)F^{''}(t)\geq\frac{5}{4}[F^{'}(t)]^2.\eqno(5.6)$$
   \\
   Note that\\
   $$[F^{-\frac{1}{4}}(t)]^{''}=-\frac{1}{4}F^{-\frac{9}{4}}(t)\left[F(t)F^{''}(t)
-\frac{5}{4}(F^{'}(t))^2\right],\eqno(5.6a)$$
\\
combining (5.5) and (5.6) together yields\\
$$[F^{-\frac{1}{4}}(t)]^{''}\leq0.$$
\\
Thus, $F^{-\frac{1}{4}}(t)$
 is concave with respect to $t$
 for sufficiently large $t$. Consequently, there exists a finite time $T^*>0$
 such that $F^{-\frac{1}{4}}(t)$ tends to 0, i.e.,
 \\
 $$\lim_{t\to T^*}F(t)=+\infty.$$
 \\
 Therefore, there exists $T<+\infty$ such that
 $$\lim_{t\to T}\left(||\phi||_{H^1(\mathbb{R^N})}^2+||\psi||_{H^1(\mathbb{R^N})}^2\right)=+\infty.
	\eqno(5.6b)$$
\\
\indent Once we prove Proposition 5.1, the proof of conclusion (1) of Theorem 2.3 will be complete. We now turn to prove Proposition 5.1 by contradiction.
\\[0.3cm]
{\bf Proof of Proposition 5.1}.~~ Assume that for all t, there holds
$$F^{'}(t)\leq0.\eqno(5.7)$$
Note that $F(t)>0$ and $F(t)$ is a convex function of $t$, one has that $F(t)$ must tend to a finite, non-negative limit $B$ as $t\to\infty$. On the other hand, since $K_1$ is an invariant set from Lemma 5.1, we can assert that $A>0$. Therefore, as $t\to\infty$,\\
$$F(t)\to A>0,\quad F^{'}(t)\to 0,\quad F^{''}(t)\to 0,$$
\\
which together with (5.5b) and $Q(\phi,\psi)<0$ yields
\\
$$\lim_{t\to\infty}\left(2a_2m_2\int_{\mathbb{R^N}}\phi_t^2dx
+a_1m_1\int_{\mathbb{R^N}}\psi_t^2dx\right)=0,\eqno(5.7a)$$\\
$$\lim_{t\to\infty}Q(\phi ,\psi ) =0.\eqno(5.7b)$$\\
According to the  conservation of energy (1.4) and (5.7a), as $t\to\infty$, there holds\\
$$	S(\phi,\psi)\leq E(\phi_{0},\psi_{0},\phi_{1},\psi_{1})<m.\eqno(5.7c)$$
\\
In view of  $Q(\phi,\psi)<0$, from Lemma 5.2 it follows that there exists $0<\alpha<1$ such that $Q(\alpha\phi,\alpha\psi)=0$. We then obtain\\
$$S(\phi,\psi)-S(\alpha\phi,\alpha\psi)\geq\frac{1}{2}Q(\phi,\psi).\eqno(5.7d)$$
\\
In view of $S(\alpha\phi,\alpha\psi)\geq m$, by (5.7b) and (5.7d), as $t\to\infty$,\\
$$S(\phi,\psi)\geq m>0,\eqno(5.7e)$$
\\
which contradicts to (5.7c). Thus there exists some $t>0$ such that $F^{'}(t)>0$, which completes the proof of Proposition 5.1.\hfill$\Box$\\
\indent Thus the proof of conclusion (1) of Theorem 2.3 is finished.\hfill$\Box$
\\[0.3cm]
\indent Next, we prove conclusion (2) of Theorem 2.3.
\\[0.3cm]
{\bf Proof of (2)}.\quad By Lemma 5.1, one has for $t\in[0,T)$,\\
$$Q(\phi,\psi)>0.\eqno(5.8)$$
\\
This together with (2.3) and (2.4) yields\\
$$S(\phi,\psi)>\frac{1}{3}\int_{\mathbb{R^N}}\left(a_2m_2|\nabla\phi|^2
+\frac{a_1m_1}{2}|\nabla\psi|^2+a_2m_2m_1^2|\phi|^2
+\frac{a_1m_1m_2^2}{2}|\psi|^2\right)dx\geq0.\eqno(5.8a)$$
\\
We further obtain\\
$$
\left.
\begin{array}{lll}
		E(\phi,\psi,\phi_t,\psi_t)&\displaystyle =S(\phi,\psi)+a_2m_2\int_{\mathbb{R^N}}\phi_t^2dx
+\frac{a_1m_1}{2}\int_{\mathbb{R^N}}\psi_t^2dx
\\[0.3cm]
&=\displaystyle E(\phi_0,\psi_0,\phi_1,\psi_1)<m.
\end{array}
\right.\eqno(5.8b)
   $$
   \\
   The above two inequalities imply that
   $$a_2m_2\int_{\mathbb{R^N}}\phi_t^2dx+\frac{a_1m_1}{2}\int_{\mathbb{R^N}}\psi_t^2dx<m,\eqno(5.8c)$$
$$S(\phi,\psi)<m.\eqno(5.8d)$$\\
On the other hand, for $t\in[0,T)$, it follow from (5.8a)-(5.8d) that\\
$$a_2m_2\int_{\mathbb{R^N}}|\nabla\phi|^2dx+\frac{a_1m_1}{2}
\int_{\mathbb{R^N}}|\nabla\psi|^2dx<m,\eqno(5.8e)$$
$$a_2m_2m_1^2\int_{\mathbb{R^N}}|\phi|^2dx+\frac{a_1m_1m_2^2}{2}\int_{\mathbb{R^N}}|\psi|^2dx<m.
 \eqno(5.8f)$$\\
That is, for $t\in[0,T)$, we establish the boundedness of $\phi,\psi$ in $H^1(\mathbb{R^N})$ and the boundedness of $\phi_t,\psi_t$ in $L^2(\mathbb{R^N})$. Consequently, we have $T=+\infty$, which implies that the solution to the Cauchy problem (1.1)-(1.2) exists globally for all $t\in[0,\infty)$. Thus, conclusion (2) of Theorem 2.3 is proved. \\
\indent Thus, the proof of Theorem 2.3 is complete. \hfill$\Box$
\section{Instability of Standing Wave}	
\indent In this section, we prove Theorem 2.4 by virtue of the results of Theorem 2.2 and Theorem 2.3.
\\[0.3cm]
{\bf Proof of Theorem 2.4}.\quad Let $(w, \xi)$ be a ground state solution of (2.1). Then it is a solution of (4.22).\\
\indent Recalling the expressions of energy and $S(u,v)$ (see (1.3),(1.4) and (2.3)),  it follows that the solution $(\phi,\psi)$ of the Cauchy problem (1.1)-(2.13) satisfy\\
$$E(\phi,\psi,\phi_t,\psi_t)=E(\phi_0,\psi_0,0,0)=S(\phi_0,\psi_0).\eqno(6.1)$$
Since $$(\lambda w,\lambda\xi)\rightarrow~(w,\xi)~~in~~H^1(\mathbb{R^N})\times H^1(\mathbb{R^N})~~\mbox{as} ~~\lambda\rightarrow 1,$$
\\
for any $\varepsilon>0$, we may choose $(\phi_{0},\psi_{0})\in H^1(\mathbb{R^N})\times H^1(\mathbb{R^N})$ with a compact support such that \\
$$\phi_0(x)=\lambda w(x),\quad \psi_0(x)=\lambda \xi(x),\quad \lambda>1,\eqno(6.2)$$
and
$$||\phi_0-w||_{H^1(\mathbb{R^N})}=(\lambda-1)||w||_{H^1(\mathbb{R^N})}<\epsilon, \quad ||\psi_0-\xi||_{H^1(\mathbb{R^N})}
=(\lambda-1)||\xi||_{H^1(\mathbb{R^N})}<\epsilon.\eqno(6.3)$$\\
Note that (6.2),  $Q(w,\xi)=0,~ \lambda>1$ and Proposition 4.4, direct calculation gives\\
$$
\left.
\begin{array}{lll}
		Q(\phi_0,\psi_0)&=Q(\lambda w,\lambda\xi) \\[0.3cm] &\displaystyle=\lambda^2\left(2a_2m_2\int_{\mathbb{R^N}}|\nabla w|^2dx+a_1m_1\int_{\mathbb{R^N}}|\nabla \xi_0|^2dx\right)
\\[0.4cm]
			&\displaystyle\quad+\lambda^2\left(2a_2m_1^2m_2\int_{\mathbb{R^N}} w_0^2dx+a_1m_2^2m_1\int_{\mathbb{R^N}}\xi^2dx\right)\\[0.4cm]			&\displaystyle\quad-6\lambda^3a_1a_2m_1m_2\int_{\mathbb{R^N}}w\xi^2dx\\[0.4cm]
			&\displaystyle<\lambda^2Q(w,\xi)=0,
\end{array}
\right.\eqno(6.4a)
   $$
 and
 $$S(\phi_0,\psi_0)=S(\lambda w,\lambda\xi)
			  <S(w,\xi)=m.	\qquad\qquad\qquad\eqno(6.4b)$$ \\
Recall (6.1), it follows that\\
  $$E(\phi_0,\psi_0,0,0)<m.	\eqno(6.5)$$ \\
  In view of conclusion (1) of Theorem 2.3, (6.4a), (6.4b) and (6.5), we complete the proof of Theorem 2.4.\hfill$\Box$
\section*{Acknowledgments}
Xiaojing Dong is partially supported by the Shandong Provincial Natural Science Foundation
(Grant No.ZR2024QA071) and the National Natural Science
Foundation of China (Grant No.12401138); Huagui Duan is partially supported by  NNSFC (No. 12271268, 12361141812) and Natural Science Foundation of Tianjin (No. 25JCZDJC01030).

\section*{Appendix A}
\indent We now prove Proposition 2.2. According to the expressions of functionals $S(u,v)$ and $Q(u,v)$~(see (2.3) and (2.4)), it suffices to show the following two propositions.
 \\[0.3cm]
 \indent Let\\
 $$\displaystyle F(u,v)\triangleq\int_{\mathbb{R}^N}vu^2 \,\mathrm{d}x.\eqno(F.1)$$
 \\
 We claim:
 \\[0.3cm]
{\bf Proposition A.I}.\quad {\sl Let $ N=5$. Then the functional $F(u,v)$ defined by (F.1) is well-defined on
	$H^1(\mathbb{R}^N)\times H^1(\mathbb{R}^N)$.}\\[0.3cm]
{\bf Proposition A.II}.\quad {\sl Let $ N=5$. Then the functional $F(u,v)$ defined by (F.1) is of class $ C^1 $ on
	$H^1(\mathbb{R}^N)\times H^1(\mathbb{R}^N)$.}
\\[0.3cm]
\indent We first prove Proposition A.I.
\\[0.3cm]
{\bf Proof of Proposition A.I}.\quad  By virtue of the Young inequality and the Gagliardo-Nirenberg inequality, we have\\
$$
\left.
\begin{array}{ll}
	\displaystyle\dfrac{N}{2}a_1a_2\int_{\mathbb{R}^N}vu^2dx&
	\displaystyle\leq C\left(\int_{\mathbb{R}^N}\dfrac{|v|^p}{p}dx
	+\int_{\mathbb{R}^N}\dfrac{|u|^{2q}}{q}dx\right)
	\\[0.5cm]
	&\displaystyle\leq C\left(\left\|v\right\|^p_{L^p(\mathbb{R}^N)}+\left\|u\right\|^{2q}_{L^{2q}(\mathbb{R}^N)}\right)
	\\[0.5cm]
	&\displaystyle\leq C\left\|u\right\|^{2q-\frac{N}{2}(2q-2)}_{L^{2}(\mathbb{R}^N)}
	\left\|\nabla u\right\|^{\frac{N}{2}(2q-2)}_{L^{2}(\mathbb{R}^N)}+C\left\|v\right\|^{p-\frac{N}{2}(p-2)}_{L^{2}(\mathbb{R}^N)}
	\left\|\nabla v\right\|^{\frac{N}{2}(p-2)}_{L^{2}(\mathbb{R}^N)},
\end{array}
\right.\eqno(A.1)
$$
 \quad\\
where\\
$$\dfrac{1}{p}+\frac{1}{q}=1,~ N=5,~ 2+\dfrac{4}{N}<p<\dfrac{2N}{N-2},~ \dfrac{2N}{N+2}<q<\frac{2N+4}{N+4}.\eqno(A.2)$$
\quad\\
On the other hand, direct computation gives \\
$$
\left\{
\begin{array}{ll}
	\displaystyle \dfrac{N}{4}\left(p-2\right)>1,\quad \dfrac{4N}{N+2}<2q<\dfrac{4N+8}{N+4},
	\\[0.5cm]
	\displaystyle \dfrac{N}{4}\left(2q-2\right)>\dfrac{N}{4}\left(\dfrac{4N}{N+2}-2\right)
	=\dfrac{N}{2}\dfrac{N-2}{N+2}>1.
\end{array}
\right.\eqno(A.3)
$$
\\
Since
$(u,v)\in H^1(\mathbb{R}^N)\times H^1(\mathbb{R}^N)$, by (A.1), we conclude the proof of Proposition A.I.\hfill$\Box$
\\[0.3cm]
\indent We then prove Proposition A.II.
\\[0.3cm]
 {\bf Proof of Proposition A.II}.\quad To do this, it suffices to verify the following two statements.
 \\[0.3cm]
 {\bf Lemma B.I}.\quad {\sl For $(u,v), (u^*,v^*)\in H^1(\mathbb{R}^N)\times H^1(\mathbb{R}^N)$,~ and for all $t_1>0,~t_2>0$, there hold:
 \\[0.3cm]
$$	 \frac{1}{t_1}\left| F(u+t_1 u^*, v) - F(u,v)-t_1\langle F_{u}(u,v), u^*\rangle \right| \to 0 ~~\mbox{as}~~   t_1\rightarrow 0  ,\eqno(A.4)
$$
and
$$
	 \frac{1}{t_2}\left| F(u, v+t_2v^*) - F(u,v)-t_2\langle F_{v}(u,v), v^*\rangle \right| \to 0 ~~\mbox{as}~~  t_2\rightarrow 0 , \eqno(A.5)
$$
\\
where
\\
$$
	F_{u}(u,v) = \int_{\mathbb{R}^N} 2vu~\mathrm{d}x, ~~ F_{v}(u,v) = \int_{\mathbb{R}^N} u^2~\mathrm{d}x,
	\eqno(A.6)
$$
$$
	\langle F_{u}(u,v), u^*\rangle = \int_{\mathbb{R}^N} 2vuu^*~\mathrm{d}x,
	\eqno(A.7)
$$
\\
$$
	\langle F_{v}(u,v),v^*\rangle = \int_{\mathbb{R}^N} u^2v^*~\mathrm{d}x.
	\eqno(A.8)
$$}
 \\[0.3cm]
 {\bf Lemma B.II}.\quad {\sl If the sequence $(u_n, v_n)$
 converges strongly to $(u,v)$
 in $H^1(\mathbb{R}^N)\times H^1(\mathbb{R}^N)$, then\\
$$
\sup_{\|u^*\|_{H^1(\mathbb{R}^N)}+\|v^*\|_{H^1(\mathbb{R}^N)}\leqslant 1}~~\left| \int_{\mathbb{R}^N} (v_n u_n -vu)u^*\,\mathrm{d}x \right| \to 0~~\mbox{as}~~  n\to \infty,
	\eqno(A.9)
$$
and
$$
  \sup_{\|u^*\|_{H^1(\mathbb{R}^N)}+\|v^*\|_{H^1(\mathbb{R}^N)}\leqslant 1}~~\left| \int_{\mathbb{R}^N} (u_n^2-u^2)v^*~\mathrm{d}x \right| \to 0 ~~\mbox{as}~~n\to\infty.\eqno(A.10)
$$}
\\[0.3cm]
\indent We then first verify Lemma B.I.
\\[0.3cm]
{\bf Proof of Lemma B.I}.\quad Let $(u,v)\in H^1(\mathbb{R}^N)\times H^1(\mathbb{R}^N), ~(u^*,v^*)\in H^1(\mathbb{R}^N)\times H^1(\mathbb{R}^N)$. Then the following estimate holds:\\
$$
\left.
\begin{array}{lll}
		&\displaystyle \left| \frac{1}{t_1} \left(F(u+t_1 u^*, v) - F(u,v)-t_1\left<F_{u}(u,v), u^*\right> \right)\right|
\\[0.4cm]
		&\displaystyle \quad=   \left| \frac{1}{t_1} \left(\int_{\mathbb{R}^N} \left[v(u+t_1u^*)^2-vu^2-2t_1vuu^*\right]\,\mathrm{d}x\right)\right|
\\[0.4cm]
		&\displaystyle \quad \leqslant  \int_{\mathbb{R}^N}\left| \frac{1}{t_1}\left[v(u+t_1u^*)^2-vu^2-2t_1vuu^*\right] \right|\,\mathrm{d}x.\nn\\
\end{array}
\right.
   $$
   \\
Note that the following estimate holds for almost every \(x\in \mathbb{R}^N\) :\\
$$
\left.
	\begin{array}{ll}	
	&\displaystyle\left| \frac{1}{t_1}\left[v(u+t_1u^*)^2-v{u}^2-2t_1vuu^*\right] \right|
\\[0.3cm]
    &\displaystyle\quad\leqslant  \left(\sup_{t\in [0,1]}\left| 2v({u+tu^*)}\right|+\left| 2(vu)\right|\right)\left| u^*\right|
    \\[0.4cm]
     &\displaystyle\quad\leqslant  C\left(\left| v \right|\left| u \right|+\left| v \right|\left| u^* \right|\right)\left| u^*\right|
     \\[0.3cm]
     &\displaystyle\quad\leqslant  C\left(\left| v \right|^\frac{N+2}{4}+\left| u \right|^\frac{N+2}{N-2}+\left| u^*\right|^\frac{N+2}{N-2}\right)\left| u^*\right|
     \\[0.3cm]
     &\displaystyle\quad\leqslant C\left(\left| u\right|^{\frac{N+2}{N-2}\cdot\frac{2N}{N+2}}+\left| u^*\right|^\frac{2N}{N-2}+\left| v\right|^{\frac{N+2}{4}\cdot\frac{8N}{N^2-4}}+\left| u^*\right|^\frac{8N}{8N-N^2+4}\right)
     \\[0.3cm]
     &\displaystyle\quad=  C\left(\left| u\right|^\frac{2N}{N-2}+\left| u^*\right|^\frac{2N}{N-2}+\left| v\right|^\frac{2N}{N-2}+\left| u^*\right|^\frac{8N}{8N-N^2+4}\right).
	\end{array}
	\right.\eqno(A.11)
$$
\\
For $N=5$, we have $2<\frac{8N}{8N-N^2+4}<\frac{2N}{N-2}$. According to the Sobolev embedding theorem,
$$
	H^1(\mathbb{R}^N)\hookrightarrow L^p(\mathbb{R}^N), \quad 2\leqslant p\leqslant \frac{2N}{N-2},
$$
one has\\
$$\left| u\right|^\frac{2N}{N-2}+\left| u^*\right|^\frac{2N}{N-2}+\left| v\right|^\frac{2N}{N-2}+\left| u^*\right|^\frac{8N}{8N-N^2+4} \in L_+^1(\mathbb{R}^N).\eqno(A.12)$$
\\
On the other hand, as $t_{1}\rightarrow 0$, the following convergence holds for almost everywhere $x\in \mathbb{R}^N$:\\
$$
	\frac{1}{t_1}\left[v({u+t_1u^*})^2-vu^2-2t_1vu{u^*}\right] \to 0~as~ t_1\to 0.
\eqno(A.13)
$$
\\
Then it follows from (A.11)-(A.12) that\\
$$
\left|\frac{1}{t_1}\left[v({u+t_1u^*})^2-vu^2-2t_1vu{u^*}\right] \right|\leq  \left| u\right|^\frac{2N}{N-2}+\left| u^*\right|^\frac{2N}{N-2}+\left| v\right|^\frac{2N}{N-2}+\left| u^*\right|^\frac{8N}{8N-N^2+4}.
\eqno(A.14)
$$
\\
Combining (A.11), (A.12) and (A.14), we obtain (A.4) by the Lebesgue Dominated Convergence Theorem.\\
\indent Moreover, direct computation yields\\
\begin{equation*}
	\begin{aligned}
		&\left| \frac{1}{t_2}\left(F(u, v+t_2v^*) - F(u,v)-t_2\left<F_{v}(u,v), {v^*}\right> \right) \right| \\[0.3cm]
		&\quad=  \left| \frac{1}{t_2} \left(\int_{\mathbb{R}^N} \left[(v+t_2v^*)u^2-vu^2-t_2u^2{v^*}\right]\,\mathrm{d}x\right)\right|
\\[0.3cm]
		&\quad= 0.
	\end{aligned}
\end{equation*}
\\
Clearly, (A.5) holds. This completes the proof of Lemma B.I.\hfill$\Box$
\\[0.3cm]
\indent Next, we prove Lemma B.II.
\\[0.3cm]
{\bf Proof of Lemma B.II}.\quad We proceed this proof in two steps.\\
\indent Let $$\Omega=\left\{x:~~|x|\leq R_{0}~~\mbox{for}~~R_{0}>0\right\}.$$
\\
\indent {\bf Step 1: ~~Estimates for $x\in \overline{\mathbb{R}^{N}\setminus \Omega}=\left\{x:~~|x|\geq R_{0}\right\}$.}
\\[0.3cm]
\indent In this case, we claim that for any
 $\varepsilon >0$, there exists $R_0>0$
 such that\\
$$\sup_{\|u^*\|_{H^1(\mathbb{R}^N)}+\|v^*\|_{H^1(\mathbb{R}^N)}\leqslant 1}~~\left|\int_{|x|\geqslant R_0} (v_n u_n -vu){u^*}\,\mathrm{d}x \right| \leqslant \varepsilon,\eqno(A.15)$$
$$\sup_{\|u^*\|_{H^1(\mathbb{R}^N)}+\|v^*\|_{H^1(\mathbb{R}^N)}\leqslant 1}~~\left| \int_{|x|\geqslant R_0} (u_n^2 -u^2){v^*}\,\mathrm{d}x \right| \leqslant \varepsilon.\eqno(A.16)$$
\\
{\renewcommand{\baselinestretch}{2.5}{\bf Proof of (A.15)}.\quad In fact, since \((u_n, v_n)\) converges strongly to $(u,v)$  in $H^1(\mathbb{R}^N)\times H^1(\mathbb{R}^N)$ ~(we may extract a subsequence if necessary), there exists $\displaystyle(u^0, v^0)\in L^{\frac{2(l+1)}{l}}(\mathbb{R}^N)\times  L^{\frac{2(l+1)}{l}}(\mathbb{R}^N)$ with $\displaystyle l+1=\frac{2N}{N-2}$ such that}
$$\left| u \right| , \left| u_n \right|\leqslant u^0~~a.e.~~ x\in \mathbb{R}^N,
	\eqno(A.17)$$
$$\left| v \right| , \left| v_n \right|\leqslant v^0~~ a.e.~~ x\in \mathbb{R}^N.
	\eqno(A.18)$$
\quad\\
Thus, for any $R>0$, by the H\"{o}lder's inequality, we have
\bea
		& \displaystyle\sup_{\|u^*\|_{H^1(\mathbb{R}^N)}+\|v^*\|_{H^1(\mathbb{R}^N)}\leqslant 1}~~\left| \int_{|x|\geqslant R} (v_n u_n -v_nu){u^*}\,\mathrm{d}x \right|\nn\\
		&\displaystyle\qquad\le \sup_{\|u^*\|_{H^1(\mathbb{R}^N)}+\|v^*\|_{H^1(\mathbb{R}^N)}\leqslant 1}~~\int_{|x|\geqslant R} \left| v_n \right| \left| u_n - u \right| \left| {u^*} \right| \,\mathrm{d}x\nn\\
		&\displaystyle\qquad\leqslant  C\sup_{\|u^*\|_{H^1(\mathbb{R}^N)}+\|v^*\|_{H^1(\mathbb{R}^N)}\leqslant 1}~~\int_{|x|\geqslant R} \left| v^0 \right| \left| u^0 \right| \left| u^*\right| \,\mathrm{d}x\nn\\
		&\displaystyle\qquad\leqslant   C\| v^0\|_{L^{\frac{2(l+1)}{l}}(|x|\geqslant R)}\| u^0\|_{L^{\frac{2(l+1)}{l}}(|x|\geqslant R)}\nn\\
&\qquad\qquad \qquad\displaystyle\times~~ \sup_{\|u^*\|_{H^1(\mathbb{R}^N)}+\|v^*\|_{H^1(\mathbb{R}^N)}\leqslant 1}\| u^*\|_{L^{\frac{2N}{N-2}}(|x|\geqslant R)}\nn\\
		&\displaystyle\qquad\leqslant   C\| v^0\|_{L^{\frac{2(l+1)}{l}}(|x|\geqslant R)} \| u^0\|_{L^{\frac{2(l+1)}{l}}(|x|\geqslant R)} .\nn
\eea
\\
Since $(u^0, v^0)\in L^{\frac{2(l+1)}{l}}(\mathbb{R}^N)\times L^{\frac{2(l+1)}{l}}(\mathbb{R}^N)$,  for any $\forall \varepsilon>0$,
 there exists $R_0>0$ such that\\
$$\sup_{\|u^*\|_{H^1(\mathbb{R}^N)}+\|v^*\|_{H^1(\mathbb{R}^N)}\leqslant 1} ~~\left| \int_{|x|\geqslant R_0} (v_n u_n -v_nu){u^*}\,\mathrm{d}x \right| \leqslant \varepsilon.\eqno(A.19)$$
\\
That is, (A.15) holds.\hfill$\sharp$
\\[0.3cm]
\indent Next, we verify that (A.16) holds.
\\[0.3cm]
{\renewcommand{\baselinestretch}{2.5} {\bf Proof of (A.16)}.\quad In fact, since $(u_n, v_n)$ converges strongly to $(u,v)$ in $H^1(\mathbb{R}^N)\times H^1(\mathbb{R}^N)$ (a subsequence may be extracted if necessary), there exists $\displaystyle(\tilde{u}^0, \tilde{v}^0)\in L^{\frac{2(l+1)}{l}}(\mathbb{R}^N)\times L^{\frac{2(l+1)}{l}}(\mathbb{R}^N)$
 such that for any $R>0$},\\
 \begin{equation*}
	\begin{aligned}
		&\sup_{\|u^*\|_{H^1(\mathbb{R}^N)}+\|v^*\|_{H^1(\mathbb{R}^N)}\leqslant 1} ~~\left| \int_{|x|\geqslant R} (u_n^2-u^2){v^*}\,\mathrm{d}x \right|
\\[0.3cm]
		 \qquad & \leqslant{} C\sup_{\|u^*\|_{H^1(\mathbb{R}^N)}+\|v^*\|_{H^1(\mathbb{R}^N)}\leqslant 1} ~~\int_{|x|\geqslant R} \left| \tilde{u}^0 \right|^2 \left| v^* \right| \,\mathrm{d}x\\[0.3cm]
		\qquad&\leqslant{}  C\| \tilde{u}^0\|_{L^{\frac{2(l+1)}{l}}(|x|\geqslant R)}^2~~ \sup_{\|u^*\|_{H^1(\mathbb{R}^N)}+\|v^*\|_{H^1(\mathbb{R}^N)}\leqslant 1}~~\| v^*\|_{L^{\frac{2N}{N-2}}(|x|\geqslant R)}\\[0.3cm]
		&\quad\leqslant{}  C\| \tilde{u}^0\|_{L^{\frac{2(l+1)}{l}}(|x|\geqslant R)}^2.
	\end{aligned}
\end{equation*}
\\
Hence, there exists  \(R_0>0\)  such that\\
$$\sup_{\|u^*\|_{H^1(\mathbb{R}^N)}+\|v^*\|_{H^1(\mathbb{R}^N)}\leqslant 1} ~~ \left|  \int_{|x|\geqslant R_0} (u_n^2 -u^2){v^*}\,\mathrm{d}x \right| \leqslant \varepsilon.\eqno(A.20)$$
\\
That is, (A.16) holds.\hfill$\sharp$
\\[0.3cm]
{\bf Step 2: Estimates for  $x\in \Omega =\{x: |x|\leqslant R_0\}$}.
\\[0.3cm]
In thsi case, we claim that if \((u_n, v_n)\) converges strongly to \((u,v)\) in
\(H^1(\Omega)\times H^1(\Omega)\), then\\
$$
 \sup_{\|u^*\|_{H^1(\Omega)}+\|v^*\|_{H^1(\Omega)}\leqslant 1}\left| \int_{\Omega} (v_n u_n -vu){u^*}\,\mathrm{d}x \right| \to 0~~\mbox{as}~~ n\to \infty,\eqno(A.21)
$$
$$\sup_{\|u^*\|_{H^1(\Omega)}+\|v^*\|_{H^1(\Omega)}\leqslant 1}\left| \int_{\Omega} (u_n^2 -u^2){v^*}\,\mathrm{d}x \right|\to 0
~~\mbox{as}~~ n\to \infty.\eqno(A.22)
$$
\\[0.3cm]
 {\bf Proof of  (A.21)}.\quad Since $(u_n, v_n)$ converges strongly to $(u,v)$~in $H^1(\Omega)\times H^1(\Omega)$ (we may take a subsequence if necessary), $(u_n, v_n)$  converges to $(u,v)$ in $\displaystyle L^{\frac{2(l+1)}{l}}(\Omega)\times L^{\frac{2(l+1)}{l}}(\Omega)$, where $ \displaystyle l+1=\frac{2N}{N-2}$. Therefore, there exists $(u^0, v^0)\in L^{\frac{2(l+1)}{l}}(\Omega)\times L^{\frac{2(l+1)}{l}}(\Omega)$~~such that\\
$$
\left| u \right| , \left| u_n \right|\leqslant u^0~~~\mbox{a.e.}~~~ x\in \Omega,\eqno(A.23)
$$
$$
\left| v \right| , \left| v_n \right|\leqslant v^0~~~\mbox{a.e.}~~~ x\in \Omega.\eqno(A.24)
$$
\\
Furthermore, by the Sobolev embedding theorem, we have
\begin{equation*}
	\begin{aligned}
		&\sup_{\|u^*\|_{H^1(\Omega)}+\|v^*\|_{H^1(\Omega)}\leqslant 1}~~\left| \int_{\Omega} (v_n u_n -vu){u^*}\,\mathrm{d}x \right|\\[0.3cm]
		&\quad\leqslant{} \sup_{\|u^*\|_{H^1(\Omega)}+\|v^*\|_{H^1(\Omega)}\leqslant 1}~~\int_{\Omega}\left(\left| v_n-v \right| \left| u \right|+ \left| u_n - u \right| \left| v_n\right|\right)\left| {u^*} \right| \,\mathrm{d}x\\[0.3cm]
		&\quad\leqslant{}  \| v_n-v\|_{L^{\frac{2(l+1)}{l}}(\Omega)}\| u^0\|_{L^{\frac{2(l+1)}{l}}(\Omega)}~~\sup_{\|u^*\|_{H^1(\Omega)}+\|v^*\|_{H^1(\Omega)}\leqslant 1}\| u^*\|_{L^{l+1}(\Omega)}\\[0.3cm]
		&\qquad+\| u_n-u\|_{L^{\frac{2(l+1)}{l}}(\Omega)}\| v^0\|_{L^{\frac{2(l+1)}{l}}(\Omega)}~~ \sup_{\|u^*\|_{H^1(\Omega)}+\|v^*\|_{H^1(\Omega)}\leqslant 1}\| u^*\|_{L^{l+1}(\Omega)}\\[0.3cm]
		&\quad\leqslant{}  C\left(\| v_n-v\|_{L^{\frac{2(l+1)}{l}}(\Omega)}+\| u_n-u\|_{L^{\frac{2(l+1)}{l}}(\Omega)}\right).
	\end{aligned}
\end{equation*}
\\
That is, (A.21) holds.\hfill$\sharp$
\\[0.3cm]
{\bf Proof of (A.22)}.\quad A similar argument yields\\
\begin{equation*}
	\begin{aligned}
		&\sup_{\|u^*\|_{H^1(\Omega)}+\|v^*\|_{H^1(\Omega)}\leqslant 1}~~\left|\int_{\Omega} (u_n^2 -u^2){v^*}\,\mathrm{d}x \right|  \\[0.3cm]
		&\quad\leqslant{} \sup_{\|u^*\|_{H^1(\Omega)}+\|v^*\|_{H^1(\Omega)}\leqslant 1}~~ \int_{\Omega} \left|u_n-u\right| \left|u_n+u\right|\left|v^*\right|\,\mathrm{d}x\\[0.3cm]
		&\quad\leqslant{} C\sup_{\|u^*\|_{H^1(\Omega)}+\|v^*\|_{H^1(\Omega)}\leqslant 1}~~\int_{\Omega} \left|u_n-u\right| \left|u^0\right|\left|v^*\right| \,\mathrm{d}x\\[0.3cm]
		&\quad\leqslant{}  C\| u_n-u\|_{L^{\frac{2(l+1)}{l}}(\Omega)}\| u^0\|_{L^{\frac{2(l+1)}{l}}(\Omega)}~~\sup_{\|u^*\|_{H^1(\Omega)}+\|v^*\|_{H^1(\Omega)}\leqslant 1}\| v^*\|_{L^{l+1}(\Omega)}\\[0.3cm]
		&\quad\leqslant{}  C\| u_n-u\|_{L^{\frac{2(l+1)}{l}}(\Omega)}.
	\end{aligned}
\end{equation*}
\\
Thus, (A.22) holds.\hfill$\sharp$
\\[0.3cm]
\indent Thus, by Step 1 and Step 2, we complete the proof of Lemma B.II. \hfill$\Box$
\\[0.3cm]
\indent Combining Lemma B.I with Lemma B.II yields the proof of Proposition A.II. Furthermore, according to \cite{BL1983}, by virtue of Proposition A.I and  Proposition A.II, we finish the proof of the Proposition 2.2. \hfill$\Box$

\end{document}